\documentclass[11pt]{article}
\usepackage{psfrag,amsmath}
\usepackage{amssymb}
\usepackage{mathrsfs}
\usepackage{amsfonts}
\usepackage{amsmath}
\usepackage{mathrsfs,color}
\usepackage{epstopdf}
\usepackage{graphicx}
\usepackage{mathrsfs,amscd,amssymb,amsthm,amsmath,url,bm}

\makeatletter

\renewcommand{\@seccntformat}[1]{{\csname the#1\endcsname}{\normalsize .}\hspace{.5em}}
\makeatother

\def \[{\begin{equation}}
\def \]{\end{equation}}

\newtheorem{thm}{Theorem}[section]
\newtheorem{prop}{Proposition}

\newtheorem{lem}[thm]{Lemma}
\newtheorem{cor}[thm]{Corollary}

\newtheorem{remark}{Remark}

\newenvironment{wst}
{\setlength{\leftmargini}{1.5\parindent}
 \begin{itemize}
 \setlength{\itemsep}{-1.1mm}}
{\end{itemize}}

\begin{document}

\setlength{\baselineskip}{15pt}
\begin{center}{\Large \bf Hermitian adjacency matrix of the second kind for mixed graphs\footnote{S.L.\, acknowledges the financial support from the National Natural Science Foundation of China (Grant Nos. 12171190, 11671164).}}
\vspace{4mm}

{\large Shuchao Li\footnote{Corresponding author. \\
\hspace*{5mm}{\it Email addresses}: lscmath@mail.ccnu.edu.cn (S.C. Li), \, ytyumath@sina.com (Y.T. Yu).},\, \, Yuantian Yu}\vspace{2mm}

Faculty of Mathematics and Statistics,  Central China Normal
University, Wuhan 430079, P.R. China
\end{center}

\noindent {\bf Abstract}:\, This contribution gives an extensive study on spectra of mixed graphs via its Hermitian adjacency matrix of the second kind { ($N$-matrix for short)} introduced by Mohar \cite{0001}. This matrix is indexed by the vertices of the mixed graph, and the entry corresponding to an arc from $u$ to $v$ is equal to the sixth root of unity $\omega=\frac{1+{\bf i}\sqrt{3}}{2}$ (and its symmetric entry is $\bar{\omega}=\frac{1-{\bf i}\sqrt{3}}{2}$); the entry corresponding to an undirected edge is equal to 1, and 0 otherwise. The main results of this paper include the following: {equivalent} conditions for a mixed graph that shares the same spectrum of its $N$-matrix with its underlying graph are given. A sharp upper bound on the spectral radius is established and the corresponding extremal mixed graphs are identified. Operations which are called two-way and three-way switchings are discussed--they give rise to some cospectral mixed graphs. We extract all the mixed graphs whose rank of its $N$-matrix is $2$ (resp. 3). Furthermore, we show that {if $M_G$ is a connected mixed graph with rank $2,$ then $M_G$ is switching equivalent to each connected mixed graph to which it is cospectral}. However, this does not hold for some connected mixed graphs with rank $3$. We identify all mixed graphs whose eigenvalues of its $N$-matrix lie in the range $(-\alpha,\, \alpha)$ for $\alpha\in\left\{\sqrt{2},\,\sqrt{3},\,2\right\}$.

\vspace{2mm} \noindent{\it Keywords:}
Mixed graph; Spectral radius; Characteristic polynomial; Switching equivalence; Cospectrality; Rank

\vspace{2mm}


\section{\normalsize Background}\setcounter{equation}{0}
Investigation on the eigenvalues of graphs has a long history. In 1965, G\"{u}nthard and Primas \cite{P4} published a paper on the spectra of trees, which probably was the first one on eigenvalues of graphs. From then on, the eigenvalue of graphs was widely used in mathematical chemistry \cite{LSG2012}, combinatorics \cite{NB1993,AEB2012,DMS1995,GGR2001,ZS2015}, code-designs theory \cite{AAH2016,CSM2017,LS} and theoretical computer science \cite{ABC2012,CDS2011} and so on. For the details, one may be referred to Guo and Mohar's contribution \cite{0005}.

In the mathematical literature, one may see that only a few works investigate the eigenvalues of directed graphs (digraphs for short). 
In the last century, the \textit{adjacency matrix} for a digraph $D$ of order $n$ was introduced, defined as an $n\times n$\, $(0, 1)$-matrix $A(D) = (a_{ij})$ with $a_{ij} = 1$ if and only if there is an arc from $v_i$ to $v_j$. This matrix attracted much attention. For the advances on this matrix, we refer the reader to the survey \cite{BR2010}. However, this matrix is not satisfied one. Clearly,  $A(D)$ is not symmetric. So many nice properties of symmetric matrices are lost for $A(D)$. A more natural definition for the adjacency matrix of a digraph was proposed by Cavers et al. \cite{MC2010}. It is called the \textit{skew-symmetric adjacency matrix} $S(D)$, in which the $(i,j)$-entry is $1$ if there is an arc from $v_i$ to $v_j$, and its symmetric entry is $-1$ (and 0 otherwise). However, this matrix works only for oriented graphs whose underlying graph is simple.

{ Reff \cite{NR} proposed the \textit{complex unit gain graph} and its \textit{adjacency matrix}. The circle group $\mathbb{T}$ is the multiplicative group of all complex numbers with norm $1$. A remarkable subgroup of $\mathbb{T}$ is the group $\mathbb{T}_n$ of the $n$-th roots of the unity. A complex unit gain graph is a triple $(G, \mathbb{T}, \varphi)$ consisting of an underlying graph $G=(V(G),E(G))$, the circle group $\mathbb{T}$ and a gain function $\varphi$, i.e., a function $\varphi$: $Dom(\varphi)\rightarrow \mathbb{T}$ such that $Dom(\varphi):= \{e_{v_iv_j}:$ \text{$e_{v_iv_j}$ belongs to $E(G)$ with ends $v_i,v_j\in $} $V(G)$\} and $\varphi(e_{v_iv_j}) =\overline{\varphi(e_{v_jv_i})}$. Here the notation $\varphi(e_{v_iv_j})$ means the gain from $v_i$ to $v_j$. In the adjacency matrix of a complex unit gain graph, the $(i,j)$-entry is the gain $\varphi(e_{v_iv_j})$ if there is an edge between $v_i$ and $v_j$, and 0 otherwise.}

Recently, Guo and Mohar \cite{0005}, and Liu and Li \cite{0008}, independently, proposed the \textit{Hermitian adjacency matrix} $H$ (or $H$-matrix for short) of a mixed graph, in which the $(i,j)$-entry is the imaginary unit ${\bf i}$ if there is an
arc from $v_i$ to $v_j$, $-{\bf i}$ if there is an arc from $v_j$ to $v_i$, 1 if $v_iv_j$ is an undirected edge, and 0 otherwise. This matrix is Hermitian and has many nice properties. Some basic theories on spectra of mixed graphs were established
via their $H$-matrices in \cite{0005,0008}. For the advances on the $H$-matrices of mixed graphs, one may be referred to \cite{SCL2020,PW2020} and in the references cited { therein}.

In 2020, Mohar \cite{0001} introduced the \textit{Hermitian adjacency matrix of the second kind} ($N$-matrix for short) of a mixed graph: each arc directed from $v_i$ to $v_j$ contributes the sixth root of unity $\omega=\frac{1 +{\bf i}\sqrt{3}}{2}$ to the $(i,j)$-entry in the matrix and contributes $\bar{\omega}=\frac{1 -{\bf i}\sqrt{3}}{2}$ to the $(j,i)$-entry; each undirected edge between $v_i$ and $v_j$ contributes 1 to the $(i,j)$- (resp. $(j,i)$-) entry, and 0 otherwise. Clearly, this novel matrix is a Hermitian matrix. It has real eigenvalues.
{One may see the $H$-matrices for mixed graphs are the adjacency matrices for complex unit gain graphs} whose gain functions take values in $\{1,\,\pm{\bf i}\}\subset\mathbb{T}_4$; whereas the $N$-matrices for mixed graphs are the adjacency matrices for complex unit gain graphs whose gain functions take values in $\{1,\,\omega,\,\bar{\omega}\}\subset\mathbb{T}_6.$ The \textit{spectrum} of a mixed graph $M_G$ is the multiset of the eigenvalues of its $N$-matrix, where the maximum modulus is called the \textit{spectral radius} of $M_G$.

{Mohar \cite{0001} showed that for a mixed bipartite graph,} its spectrum is symmetric about 0; he established some relationships between the spectral radius and the largest eigenvalue of this new matrix. In this article, we investigate some basic properties of the $N$-matrix, which may be viewed as a continuance of Mohar's work \cite{0001}.
A sharp upper bound on the spectral radius of mixed graphs is established and the corresponding extremal mixed graphs are identified (see Section 3).

{Two mixed graphs are \textit{cospectral} if they have the same spectrum. We mainly consider the cospectrality between two mixed graphs which have the same underlying graph.} Operations which are called two-way and three-way switchings are discussed--they give rise to a large number of cospectral mixed graphs. {However, there are cases in which these operations yield very few (sometimes as few as zero) switching equivalent, non-isomorphic mixed graphs.} Some equivalent conditions for a mixed graph that shares the same spectrum with its underlying graph are deduced (see Section 4).

It is interesting to study the rank of {the $N$-matrix for a mixed graph}. We extract all the mixed graphs whose rank equals $2$ (resp. $3$). Furthermore, we show that {if $M_G$ is a connected mixed graph with rank $2,$ then $M_G$ is switching equivalent to each connected mixed graph to which it is cospectral}. However, this does not hold for some connected mixed graphs with rank $3$. {These kinds of questions are located in Section 5.}

Despite many {unusual} properties that the $N$-matrix exhibits,
{it is challenging to derive combinatorial structure of a mixed graph from its eigenvalues.}
In Section 6, we find all mixed graphs whose eigenvalues lie in the range $(-\alpha, \alpha)$ for
$\alpha\in\left\{\sqrt{2},\,\sqrt{3},\,2\right\}$. 

\section{\normalsize Some definitions and preliminaries}

In this paper, we consider only simple and finite graphs. For graph theoretic notation and terminology not defined here, we refer to \cite{DW1996}.

Let $G =(V(G), E(G))$ be a graph with vertex set $V(G)$ and edge set $E(G)$. {The number of vertices $n=|V(G)|$ in a graph is called the \textit{order}.} We say that two vertices $i$ and $j$ are \textit{adjacent} (or \textit{neighbours}) if they are joined by an edge and we write it as $i\sim j$. The {\textit{degree} $d_G(u)$} of a vertex $u$ (in a graph $G$) is the number of edges incident with it. In particular, {the \textit{maximum degree} is} denoted by $\Delta(G)$. The \textit{set of neighbours} of a vertex $u$ is denoted by $N_G(u)$. A $k$-\textit{partite graph} is a graph whose set of vertices is decomposed into $k$ disjoint sets such that no two vertices within the same set are adjacent. A \textit{complete} $k$-\textit{partite graph} is a $k$-partite graph in which two vertices are adjacent if and only if they belong to different sets. As usual, let $P_n, C_n$ and $K_n$ denote the path, cycle and complete graph on $n$ vertices, respectively. { And let $K_{n_1,n_2,\ldots,n_k}$ denote the complete $k$-partite graph with the {orders} of partite sets being $n_1,\,n_2,\,\ldots,\,n_k,$ respectively. The \textit{girth} of a graph is the length of the shortest cycle contained in it.} We use $kG$ to denote the disjoint union of $k$ copies of $G$.

A \textit{mixed graph} $M_G$ is obtained from a simple graph $G$, the underlying graph of $M_G$, by orienting each edge of some subset $E_0\subseteq E(G)$. 
Thus, mixed graphs are the generalizations of simple graphs and directed graphs. A mixed graph $M_{G'}$ is a mixed subgraph of $M_G$ if $G'$ is a subgraph of $G$ and the direction of each edge in $M_{G'}$ coincides with that in $M_G$. {For a vertex subset $V'$ of $V(G),\, M_G[V']$ is a mixed subgraph of $M_G$ induced on $V'.$ A mixed graph $M_{G'}$ is called an induced subgraph of $M_G$ if there is a vertex subset $V'$ of $V(G)$ such that $M_{G'}\cong M_G[V'].$ The \textit{order} 
of $M_G$ is exactly the order 
of $G$. A mixed graph is called to be \textit{connected} if its underlying graph is connected.

We write an undirected edge as $\{u,v\}$ and a directed edge (or an arc) from $u$ to $v$ as $\overrightarrow{uv}.$ Usually, we denote an edge of $M_G$ by $uv$ if we {are not concerned} whether it is directed or not. Then $M_G-u,\,M_G-uv$ are the mixed graphs obtained from $M_G$ by deleting the vertex $u \in V(G)$ and the edge $uv \in E(M_G)$, respectively. This notation is naturally extended if more than one vertex or edge are deleted.

Given a mixed graph $M_G=(V(M_G), E(M_G)),$ let $N_{M_G}^0(v)=\{u\in V(M_G):\, \{u,v\}\in E(M_G)\}$, $N_{M_G}^+(v)=\{u\in V(M_G):\,\overrightarrow{vu}\in E(M_G)\}$ and $N_{M_G}^-(v)=\{u\in V(M_G):\,\overrightarrow{uv}\in E(M_G)\}$. Clearly, $N_G(v)=N_{M_G}^0(v)\cup N_{M_G}^+(v)\cup N_{M_G}^-(v)$. In our context, two vertices $u,v$ in a mixed graph are called to be adjacent if they are adjacent in its underlying graph and we also denote it by $u\sim v.$ The degree of a vertex in a mixed graph $M_G$ is defined to be the degree of this vertex in the underlying graph $G$.

The \textit{Hermitian adjacency matrix of the second kind}, written as {$N(M_G)=(n_{st})$}, of a mixed graph $M_G$ was proposed by Mohar \cite{0001}. It is defined as
$$
{n_{st}}=\left\{
         \begin{array}{cl}
           \omega, & \text{if $\overrightarrow{u_su_t}$ is an arc from $u_s$ to $u_t$;} \\[5pt]
    \bar{\omega}, & \text{if $\overrightarrow{u_tu_s}$ is an arc from $u_t$ to $u_s$;} \\[5pt]
    1, & \text{if $\{u_s,u_t\}$ is an undirected edge;} \\[5pt]
    0, & \text{otherwise},
         \end{array}
       \right.
$$
{where $\omega=\frac{1+ {\bf i}\sqrt{3}}{2}$ is the sixth root of unity, $\bar{\omega}$ is the complex conjugate of $\omega$. For convenience, in this paper we abbreviate this novel matrix as $N$-matrix. {With this matrix formulation, also $N^c = J -N -I$, where $N^c$ is the $N$-matrix
of the complement, $J$ is the all one matrix and $I$ is the identity matrix, the usual relation between the eigenvalues of a graph and
those of its complement carries over.}

The main reason why the sixth root of unity is natural is that $\omega+\bar{\omega}=1,$ {and so an undirected edge can be seen as two arcs with opposite directions.} In considering the relationship between the order of a mixed graph and the multiplicities of its ``weighted" Hermitian adjacency eigenvalues, the ``weights" $-\omega$ and $-\bar{\omega}$ also play an important role; see \cite{ALG}. Furthermore,} the sixth root of unity emerges realistically across applications. It appears in the definition of Eisenstein integers; in relation to matroid theory, the sixth root matroids play a special role next to regular and binary matroids; see \cite{0004,0006} for details. 

{The \textit{rank} of $M_G$ is the rank of $N(M_G)$. The \textit{characteristic polynomial} of $N(M_G), P_{M_G}(x) = \det(xI -N(M_G))$, is also called the characteristic polynomial of $M_G$, while its roots are just the eigenvalues of $M_G$.}

Note that $N(M_G)$ is Hermitian, that is, $N^*(M_G) =N(M_G)$, where $N^*(M_G)$ denotes the conjugate transpose of $N(M_G)$. Then its eigenvalues are real. The collection of eigenvalues of $M_G$ (with repetition) is called the \textit{spectrum} of $M_G$. We denote the eigenvalues of $M_G$ by
${ \lambda_1 { \geq} \lambda_2 { \geq} \cdots { \geq} \lambda_n,}
$
where $n$ is the order of $M_G$. Two mixed graphs are called \textit{cospectral} if they have the same spectrum. The \textit{spectral radius} of ${M_G}$, written as $\rho({ M_G})$, is defined as
$$
\rho({M_G})=\max\{|\lambda_1|, |\lambda_n|\}.
$$

A \textit{mixed cycle} is a mixed graph whose underlying graph is a cycle. A mixed cycle is \textit{even} (resp. \textit{odd}) if its order is even (resp. odd). Let ${M_G}$ be a mixed graph, and let $M_C=v_1v_2v_3\cdots v_{l-1}v_lv_1$ be a mixed cycle of ${M_G}$. {Note that $M_C$ cannot have repeated vertices}. Then the weight of $M_C$ in a direction is defined by
$$
wt(M_C)=n_{12}n_{23}\cdots n_{(l-1)l}n_{l1},
$$
where $n_{jk}$ is the $(v_j,v_k)$-entry of ${N(M_G)}$. It is easy to obtain $wt(M_C)\in \mathbb{T}_6.$ Note that if, for one direction, the weight of a mixed cycle is $\alpha$, then for the reversed direction its weight is $\bar{\alpha}$, the conjugate of $\alpha$. For a mixed cycle $M_C$, it is \textit{positive} (resp. \textit{negative}) if $wt(M_C)=1$ (resp. $-1$); it is \textit{semi-positive} if $wt(M_C)\in \{\omega,\bar{\omega}\}$, whereas  it is \textit{semi-negative} if $wt(M_C)\in \{-\omega,-\bar{\omega}\}$. An example of positive (resp. semi-positive, semi-negative, and negative) mixed cycle is depicted in Figure \ref{fig1}. Furthermore, we call a mixed graph ${M_G}$ \textit{positive}, if each mixed cycle of ${M_G}$ is positive.
\begin{figure}[h!]\label{fig1}
\begin{center}
\includegraphics[width=80mm]{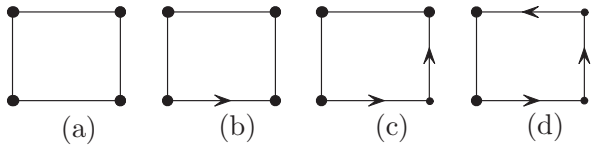} \\
  \caption{The mixed cycles: (a) is positive; (b) is semi-positive; (c) is semi-negative and (d) is negative.}\label{fig1}
\end{center}
\end{figure}

Further on we need the following preliminary results. {For a complex number $c$, let $\mathfrak{R}(c)$ denote the real part of $c.$}

{ Reff \cite[Lemma 2.2]{NR2016} characterized a sufficient condition for the cospectrality of complex unit gain graphs. As a corollary, we have the following result.
\begin{cor}\label{cor2.1}
Let ${M_G}$ and ${M'_G}$ be mixed graphs with the same underlying graph $G.$ If for every cycle $C$ in $G$, $\mathfrak{R}(wt(M_C))=\mathfrak{R}(wt(M'_C)),$ then ${M_G}$ and ${M'_G}$ are cospectral, where $M_C$ (resp. $M'_C$) is a mixed cycle in $M_G$ (resp. $M'_G$) whose underlying graph is $C.$
\end{cor} }
Let $M_G$ be a mixed graph with connected components $M_{G_1},\, M_{G_2},\, \ldots,\, M_{G_t}.$ Then $N(M_G)$ can be written as
$$
N(M_G)=\left(
                   \begin{array}{ccc}
                     N(M_{G_1}) & \  & \  \\
                     \  & \ddots & \  \\
                     \  & \  & N(M_{G_t}) \\
                   \end{array}
                 \right).
$$
Hence the following result is clear.
\begin{lem}\label{lem2.2}
Let $M_G$ be a mixed graph with connected components $M_{G_1},\, M_{G_2},\, \ldots,\, M_{G_t}.$ Then
$$
P_{M_G}(x)=\prod_{j=1}^tP_{M_{G_j}}(x).
$$
\end{lem}
{The following result establishes a relationship between the eigenvalues and the number of edges of a mixed graph.}
\begin{lem}\label{lem2.3}
Let $M_G$ be an $n$-vertex mixed graph with $m$ edges and let $\lambda_1,\,\lambda_2,\,\cdots,\,\lambda_n$ be its eigenvalues. Then
$
\sum_{j=1}^n\lambda_j^2=2m.
$
\end{lem}
\begin{proof}
Let $N=N(M_G)$. Since $N$ is Hermitian and has only entries $0,\, 1,\, \omega$ and $\bar{\omega}$, we have
$$
N_{uv}N_{vu}=N_{uv}\overline{N_{uv}}=1
$$
whenever $N_{uv}\neq0$. This implies that the $(u,u)$ entry in $N^2$ is the degree of $u$ in $G$. Hence
$$
\sum_{j=1}^n\lambda_j^2=tr(N^2)=\sum_{u\in V(G)}(N^2)_{uu}=\sum_{u\in V(G)}d_G(u)=2m,
$$
where $tr(N^2)$ denotes the trace of the matrix $N^2.$
\end{proof}
Suppose that $a_1\geq a_2\geq\cdots\geq a_n$ and $b_1\geq b_2\geq\cdots\geq b_{n-t}$ (where $t\geq1$ is an integer) be two sequences of real numbers. We say that the sequences $a_l\, (1\leq l\leq n)$ and $b_j\, (1\leq j\leq n-t)$ \textit{interlace} if for every $s=1,\ldots,n-t$, we have
$
a_s\geq b_s\geq a_{s+t}.
$
The following interlacing theorem is well-known; { see \cite[Theorem~4.1]{0005}.}
\begin{thm}\label{thm2.4}
Let $B$ be a Hermitian matrix and $B'$ be its principal submatrix. Then the eigenvalues of $B'$ interlace those of $B$.
\end{thm}
Theorem~\ref{thm2.4} implies that the eigenvalues of any induced mixed subgraph interlace those of the mixed graph itself.
\begin{cor}\label{cor2.5}
The eigenvalues of an induced mixed subgraph interlace the eigenvalues of the mixed graph.
\end{cor}

An \textit{elementary mixed graph} is a mixed graph such that every component is either an (oriented) edge or a mixed cycle. The \textit{rank} of a simple graph $G$ is defined by $r(G)=n-t(G)$, where $n$ and $t(G)$ are the order and number of components of $G$, respectively.

The coefficients of the characteristic polynomial of the adjacency matrix for a complex unit gain graph have been determined; see \cite[Theorem 2.7]{AS2019}.
As a corollary, given an $n$-vertex mixed graph $M_G$, we can obtain a description of all the coefficients of the characteristic polynomial $P_{M_G}(x)$.} Let
\[\label{eq:3.01}
P_{M_G}(x)=x^n+c_1x^{n-1}+c_2x^{n-2}+\cdots+c_{n-1}x+c_n,
\]
where $c_1,\ldots, c_n$ are real.
\begin{thm}\label{thm3.2}
Let $M_G$ be a mixed graph of order $n$. Then the coefficients of the characteristic polynomial $P_{M_G}(x)$ in \eqref{eq:3.01} {are} given by
$$
  c_k=\sum_{M_{G'}}(-1)^{-k+r(G')+l_{sn}(M_{G'})+l_n(M_{G'})}\cdot 2^{l_p(M_{G'})+l_n(M_{G'})},
$$
where the summation is over all elementary mixed subgraphs $M_{G'}$ of $M_G$ with $k$ vertices and underlying graph $G',\,l_p(M_{G'}),\,l_n(M_{G'}),\,l_{sn}(M_{G'})$ are the number of positive, negative, semi-negative cycles in $M_{G'}$, respectively.
\end{thm}
{ \begin{proof}
According to \cite[Theorem 2.7]{AS2019}, it is sufficient to show that for every $k$-vertex elementary mixed subgraph $M_{G'}$ of $M_G,$
$$
(-1)^{t(M_{G'})}2^{l(M_{G'})}\prod_{M_C\in \mathcal{C}(M_{G'})}\mathfrak{R}(wt(M_C))=(-1)^{-k+r(G')+l_{sn}(M_{G'})+l_n(M_{G'})}\cdot 2^{l_p(M_{G'})+l_n(M_{G'})},
$$
where $t(M_{G'})$ is the number of components of $M_{G'},\,l(M_{G'})$ is the number of mixed cycles in $M_{G'},\, \mathcal{C}(M_{G'})$ is the collection of all mixed cycles in $M_{G'}.$ In fact, according to the definitions of positive, negative, semi-positive and semi-negative cycles, one has
\begin{align*}
&(-1)^{t(M_{G'})}2^{l(M_{G'})}\prod_{M_C\in \mathcal{C}(M_{G'})}\mathfrak{R}(wt(M_C)) \\
=&(-1)^{t(M_{G'})}2^{l(M_{G'})}\cdot1^{l_p(M_{G'})}\cdot(-1)^{l_n(M_{G'})}\cdot(\frac{1}{2})^{l_{sp}(M_{G'})}\cdot(-\frac{1}{2})^{l_{sn}(M_{G'})} \\
=&(-1)^{k-r(G')+l_{sn}(M_{G'})+l_n(M_{G'})}\cdot 2^{l_p(M_{G'})+l_n(M_{G'})},
\end{align*}
where $l_{sp}(M_{G'})$ is the number of semi-positive cycles in $M_{G'}$.
\end{proof}
}
From Theorem~\ref{thm3.2} we can deduce that $c_1=0$ and $c_2=-|E(M_G)|$ for each mixed graph $M_G$. As $M_G$ has no elementary subgraph of order $1$, and has $|E(M_G)|$ elementary subgraphs of order $2$, each of which is an edge and hence contributes $-1$ to $c_2$. In \cite{0001}, Mohar showed that if $M_G$ is a mixed graph whose underlying graph is bipartite, then the spectrum of $M_G$ is symmetric about $0$. This can be easily seen from Theorem~\ref{thm3.2}, as $M_G$ has no elementary subgraph of odd order, $c_k=0$ if $k$ is odd.

In the following, we will give two recurrence relations for $P_{M_G}(x)$, which are similar to those of adjacency matrices of simple graphs \cite[Section 2]{0002} and those of Hermitian adjacency matrices for mixed graphs \cite{0016}.
\begin{thm}\label{thm3.6}
Let $M_G$ be a mixed graph, and let $u$ be a vertex of $M_G$. Then
\[\label{eq:3.1}
 P_{M_G}(x)=xP_{M_G-u}(x)-\sum_{v\sim u}P_{M_G-v-u}(x)-\sum_{Z\in \mathscr{C}(u)}\left(wt(Z)+\overline{wt(Z)}\right)P_{M_G-V(Z)}(x),
\]
where $\mathscr{C}(u)$ is the set of mixed cycles containing $u$, $wt(Z)$ is the weight of $Z$ in a direction.
\end{thm}
\begin{proof}
The proof follows the same line as the proof of \cite[Theorem 2.3.4]{0002}, the key difference is that when $u$ is contained in a mixed cycle of an elementary subgraph of $M_G,$ we need to discuss the weight of this mixed cycle.


For an elementary subgraph $M_{G'}$ of $M_G$ on $k$ vertices, if $u$ lies in a mixed cycle $Z$ of $M_{G'}$, then take $M_{G''}=M_{G'}-V(Z),$ regarded as an elementary subgraph of $M_G-V(Z)$.

By applying Theorem~\ref{thm3.2}, we can show that if $M_{G'}$ contributes $c$ to the coefficient of $x^{n-k}$ on the left of equation \eqref{eq:3.1}, i.e., $(-1)^{-k+r(G')+l_{sn}(M_{G'})+l_n(M_{G'})}\cdot 2^{l_p(M_{G'})+l_n(M_{G'})}=c,$ then $M_{G''}$ contributes $c$ to the coefficient of $x^{n-k}$ on the right of equation \eqref{eq:3.1}.

In fact, if $|V(Z)|=s$, then the contribution of $M_{G''}$ to the coefficient of $x^{(n-s)-(k-s)}\, (=x^{n-k})$ in $P_{M_G-V(Z)}(x)$ is
$$
(-1)^{k-s}\cdot(-1)^{r({G''})+l_{sn}(M_{G''})+l_n(M_{G''})}\cdot 2^{l_p(M_{G''})+l_n(M_{G''})}.
$$
Note that,
$$
r({G'})-r({G''})=(|V({G'})|-t({G'}))-(|V({G''})|-t({G''}))=s-1.
$$
If $Z$ is a positive cycle, i.e., $wt(Z)=1,$ then $l_p(M_{G'})-l_p(M_{G''})=1$ and $l_n(M_{G'})=l_n(M_{G''}),\,l_{sn}(M_{G'})=l_{sn}(M_{G''}).$ Hence,
\begin{align*}
&(-1)^{k-s}\cdot(-1)^{r({G''})+l_{sn}(M_{G''})+l_n(M_{G''})}\cdot 2^{l_p(M_{G''})+l_n(M_{G''})}\\
=&-(-1)^k\cdot\frac{1}{2}(-1)^{r({G'})+l_{sn}(M_{G'})+l_n(M_{G'})}\cdot 2^{l_p(M_{G'})+l_n(M_{G'})}\\
=&-\frac{1}{2}c.
\end{align*}
This gives that the contribution of $M_{G''}$ to the coefficient of $x^{n-k}$ in $\left(wt(Z)+\overline{wt(Z)}\right)P_{M_G-V(Z)}(x)$ is $-c$.
Similarly, we can prove that if $Z$ is a negative, semi-positive or semi-negative cycle, then the contribution of $M_{G''}$ to the coefficient of $x^{n-k}$ in $\left(wt(Z)+\overline{wt(Z)}\right)P_{M_G-V(Z)}(x)$ is also $-c$. Besides, $M_{G''}$ does not contribute to the coefficient of $x^{n-k}$ in the remaining terms on the right of equation \eqref{eq:3.1}, hence $M_{G''}$ contributes $c$ to the coefficient of $x^{n-k}$ on the right of equation \eqref{eq:3.1}.
\end{proof}
\begin{thm}\label{thm3.7}
Let $M_G$ be a mixed graph, and let $uv$ be a mixed edge of $M_G$. Then
$$
P_{M_G}(x)=P_{M_G-uv}(x)-P_{M_G-v-u}(x)-\sum_{Z\in \mathscr{C}(uv)}\left(wt(Z)+\overline{wt(Z)}\right)P_{M_G-V(Z)}(x),
$$
where $\mathscr{C}(uv)$ is the set of mixed cycles containing $uv$, $wt(Z)$ is the weight of $Z$ in a direction.
\end{thm}
\begin{proof}
The proof is similar to the proof of Theorem~\ref{thm3.6}, and we omit it here.
\end{proof}
\begin{cor}\label{cor3.8}
Let $M_G$ be a mixed graph with $uv$ being a cut edge of its underlying graph, and let $M_{G_1},\,M_{G_2}$ be two components of $M_G-uv$ with $u\in V(M_{G_1}),\,v\in V(M_{G_2}).$ Then
$$
P_{M_G}(x)=P_{M_{G_1}}(x)P_{M_{G_2}}(x)-P_{M_{G_1}-u}(x)P_{M_{G_2}-v}(x).
$$
\end{cor}
\begin{proof}
According to Theorem~\ref{thm3.7}, we have
$$
P_{M_G}(x)=P_{M_G-uv}(x)-P_{M_G-v-u}(x),
$$
as $uv$ is contained in no mixed cycle of $M(G)$. By Lemma~\ref{lem2.2},
$$
P_{M_G-uv}(x)=P_{M_{G_1}}(x)P_{M_{G_2}}(x),\ \ P_{M_G-v-u}(x)=P_{M_{G_1}-u}(x)P_{M_{G_2}-v}(x).
$$
This completes the proof.
\end{proof}
This result is the same as the corresponding result for simple graphs which has been proved in \cite[Section 2]{0002} by another method. More recurrence relations for $P_{M_G}(x)$ which are the same as the case of simple graphs can be seen in \cite[Section 2]{0002}.

\section{\normalsize An upper bound for the spectral radius}\setcounter{equation}{0}
In this section, we show that $\rho(M_G)$ is bounded above by $\Delta(G)$ and when $G$ is connected, we characterize the mixed graphs attaining this bound.
\begin{thm}\label{thm4.1}
Let $M_G$ be an $n$-vertex mixed graph whose underlying graph is $G$. Then $\rho(M_G)\leq\Delta(G)$. When $G$ is connected, the equality holds if and only if $G$ is $\Delta(G)$-regular and one can partition $V(M_G)$ into six (possibly empty) parts { $V_1,V_{-1},V_{\omega}, V_{\bar{\omega}},V_{-\omega},V_{-\bar{\omega}}$} such that one of the followings holds:
\begin{wst}
\item[{\rm (i)}]The induced mixed graph $M_G[V_j]$ contains only undirected edges for all $j\in\mathbb{T}_6$ and each of the rest edges in ${ E(M_G)\setminus (\bigcup_{j\in\mathbb{T}_6}E(M_G[V_j]))}$ is an arc $\overrightarrow{uv}$ satisfying $u\in V_j$ and $v\in V_{\bar{\omega}\cdot j}$ for some $j\in\mathbb{T}_6$; see Figure~$\ref{fig2}$.
\item[{\rm (ii)}]The induced mixed graph $M_G[V_j]$ is an independent set for all $j\in\mathbb{T}_6$ ; every undirected edge $\{u,v\}$ of $M_G$ satisfies $u\in V_j$ and $v\in V_{-j}$ for some $j\in\mathbb{T}_6$, and every arc $\overrightarrow{uv}$ of $M_G$ satisfies $u\in V_j$ and $v\in V_{-\bar{\omega}\cdot j}$ for some $j\in\mathbb{T}_6$; see Figure~$\ref{fig2}$.
\end{wst}
\begin{figure}[h!]
\begin{center}
\includegraphics[width=84mm]{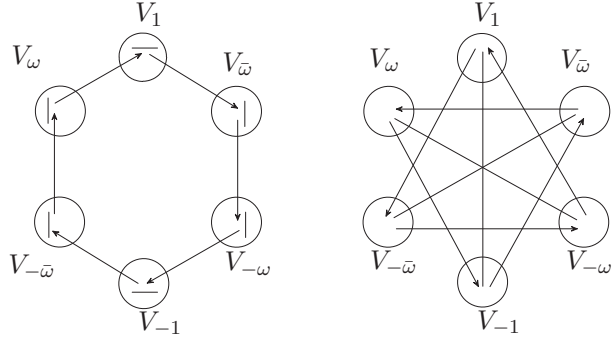} \\
  \caption{Items (i) and (ii) of Theorem~\ref{thm4.1}.}\label{fig2}
\end{center}
\end{figure}
\end{thm}
\begin{proof}
Let $N=N(M_G)$ and let $\mathbf{x}=(x_1, x_2, \ldots, x_n)^T$ be an eigenvector corresponding to the eigenvalue $\lambda$ of $N$. Associate a labeling of vertices of $M_G$ (with respect to $\mathbf{x}$) in which $x_i$ is a label of $v_i$. Without loss of generality, let $|x_1|=\max\{|x_i|: 1\leq i\leq n\}.$ On the one hand,  we consider the first entry of $N\mathbf{x}$:
$$
(N\mathbf{x})_1=\left(\sum_{v_i\in N_{M_G}^0(v_1)}x_i\right)+\omega\left(\sum_{v_j\in N_{M_G}^+(v_1)}x_j\right)+\bar{\omega}\left(\sum_{v_k\in N_{M_G}^-(v_1)}x_k\right).
$$
On the other hand, from $N\mathbf{x}=\lambda \mathbf{x}$, we obtain
\[\label{eq:4.1}
  (N\mathbf{x})_1=\lambda x_1.
\]
This gives us
\begin{eqnarray}
|\lambda x_1|&=&|(N\mathbf{x})_1|\notag\\
&=&\left|\left(\sum_{v_i\in N_{M_G}^0(v_1)}x_i\right)+\omega\left(\sum_{v_j\in N_{M_G}^+(v_1)}x_j\right)+\bar{\omega}\left(\sum_{v_k\in N_{M_G}^-(v_1)}x_k\right)\right| \notag
\end{eqnarray}
\begin{eqnarray}
&\leq &\left(\sum_{v_i\in N_{M_G}^0(v_1)}|x_i|\right)+|\omega|\left(\sum_{v_j\in N_{M_G}^+(v_1)}|x_j|\right)+|\bar{\omega}|\left(\sum_{v_k\in N_{M_G}^-(v_1)}|x_k|\right) \label{eq:4.2}\\
&\leq&d_G(v_1)|x_1| \label{eq:4.3}\\ \label{eq:4.4}
&\leq&\Delta(G)|x_1|.
\end{eqnarray}
Hence, $|\lambda|\leq\Delta(G)$. Note that $\lambda$ is an arbitrary eigenvalue of $M_G$, and by the definition of spectral radius, we have  $\rho(M_G)\leq\Delta(G)$.

In what follows, we characterize all the mixed graphs attaining this bound if the underlying graph $G$ is connected.

Note that $\rho(M_G)=\Delta(G)$ holds if and only if equalities above must hold throughout. 
We see that the equality in \eqref{eq:4.4} holds if and only if $d_G(v_1)=\Delta(G)$, whereas the equality in \eqref{eq:4.3} holds if and only if
\[
  \text{$|x_k|=|x_1|$\, \, for all\, \, $v_k\in N_G(v_1)$.}
\]
Since the choice of $v_1$ is arbitrary among all vertices attaining the maximum absolute value in $\mathbf{x}$, we may apply the same discussion to any vertex adjacent to $v_1$ in $G$. Note that $G$ is connected. Therefore, $G$ is $\Delta(G)$-regular
and $|x_k|=|x_1|$ for all $v_k\in V(M_G)$.

We may normalize $\mathbf{x}$ such that $x_1=1$. Hence, $|x_i|=1$ for $i\in\{1,2,\ldots,n\}$. The inequality in \eqref{eq:4.2} follows from the triangle inequality for sums of complex numbers, and so the equality holds if and only if every complex number in the following set $S$ has the same argument, where
\[\label{eq:4.6}
 S=\left\{x_i: v_i\in N_{M_G}^0(v_1)\right\}\cup \left\{\omega\cdot x_j: v_j\in N_{M_G}^+(v_1)\right\}\cup \left\{\bar{\omega}\cdot x_k: v_k\in N_{M_G}^-(v_1)\right\}.
\]
Together with \eqref{eq:4.1}, the equality in \eqref{eq:4.2} holds if and only if every complex number in $S$ has the same argument as $\lambda x_1$. There are three cases for $\lambda$: $\lambda=0,\, \lambda>0$ or $\lambda<0$.  Since $\rho(M_G)\le \Delta(G)$, and the only mixed graph with $\rho(M_G)=0$ is the empty graph, it suffices to consider the following two cases.

{\bf Case 1.} $\lambda>0.$ In this case, if $\rho(M_G)=\Delta(G)$, then together with \eqref{eq:4.1} we have $(N\mathbf{x})_1=\Delta(G)x_1$. Combining with \eqref{eq:4.6} we deduce that every complex number in $S$ is just $x_1$ and is thus equal to $1$. We conclude that
\begin{align}\notag
x_i=\left\{
  \begin{array}{cl}
    1, & \text{if}\, v_i\in N_{M_G}^0(v_1); \\[5pt]
    \bar{\omega}, & \text{if}\, v_i\in N_{M_G}^+(v_1);\\[5pt]
    \omega, & \text{if}\, v_i\in N_{M_G}^-(v_1).
  \end{array}
\right.
\end{align}
Similar argument can be applied to { each $x_j\in \mathbb{T}_6\backslash\{1\}$.} From this we conclude that $V(M_G)$ is partitioned into
$$
V_1\cup V_{-1}\cup V_{\omega}\cup V_{\bar{\omega}}\cup V_{-\omega}\cup V_{-\bar{\omega}}
$$
according to the value of $x_j$, and so item (i) holds.

{\bf Case 2.} $\lambda<0.$ In this case, if $\rho(M_G)=\Delta(G)$, then together with \eqref{eq:4.1} we have $(N\mathbf{x})_1=-\Delta(G)x_1$. Combining with \eqref{eq:4.6} we obtain that every complex number in $S$ is just $-x_1$ and thus equals $-1$. { With the same discussion as that of Case 1, we conclude that, in this case,} $V(M_G)$ has a partition
$$
V_1\cup V_{-1}\cup V_{\omega}\cup V_{\bar{\omega}}\cup V_{-\omega}\cup V_{-\bar{\omega}}
$$
satisfying item (ii).

Now, we consider the converse for the two cases of the theorem. Let $M_G$ be a mixed graph whose underlying graph is $k$-regular. Assume that $V(M_G)$ has a partition $\bigcup_{j\in \mathbb{T}_6}V_j$ satisfying item (i) or (ii).

Let $\mathbf{x}$ be the vector indexed by the vertices of $M_G$ such that $x_i=j$ if $v_i\in V_j$, where $j\in \mathbb{T}_6$. Then it is easy to see that for every vertex $v_i$ we have $(N\mathbf{x})_i=kx_i$ (by item (i)) or $(N\mathbf{x})_i=-kx_i$ (by item (ii)). Thus $\mathbf{x}$ is an eigenvector of $N$ corresponding to the eigenvalue $k$ or $-k$, and so $\rho(M_G)=k=\Delta(G)$. Then the bound is tight as claimed.
\end{proof}
\begin{figure}[h!]\label{fig3}
\begin{center}
\includegraphics[width=25mm]{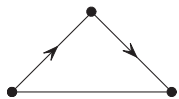} \\
  \caption{The semi-negative triangle.}
\end{center}
\end{figure}
For simple graphs, $\rho(G)$ is always larger or equal to the average degree. However, for mixed graphs, $\rho(M_G)$ can be smaller than the minimum degree of $G$. An example is the semi-negative triangle shown in Figure~$3$, whose characteristic polynomial is $x^3-3x+1$ (based on  Theorem~\ref{thm3.2}). Clearly, its spectral radius is less than $2$, while the minimum degree of its underlying graph is $2$. Of course, this anomaly is also confirmed by Theorem~\ref{thm4.1}, since the semi-negative triangle shown in Figure~\ref{fig3} does not have the structure as depicted in Figure~\ref{fig2}. {This phenomenon is similar to what happens w.r.t. $H$-matrix of a mixed graph, which has been mentioned in \cite{0005}.}

\section{\normalsize Switching equivalence and cospectrality}\setcounter{equation}{0}

In this section, we focus on properties of mixed graphs that are cospectral by introducing some
operations on mixed graphs that preserve the spectrum. In particular, we are inspired to
study mixed graph operations that preserve the spectrum and conserve the underlying
graph. We try to demonstrate the spectral information about the underlying graph by looking
at some spectrum preserving operations that do not change the underlying graph.

\begin{figure}[h!]\label{fig05}
\begin{center}
\includegraphics[width=90mm]{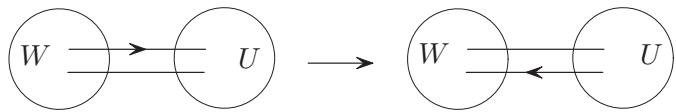} \\
  \caption{A two-way switching.}
\end{center}
\end{figure}

{ Let $W\subset V(M_G)$ be non-empty, and let $U:=V(M_G)\backslash W.$ If $W$ is such that $M_G$ contains no arc $\overrightarrow{uw}$ with $u\in U$ and $w\in W,$ then a \textit{two-way switching} is said to be the operation that replaces every arc $\overrightarrow{wu}\,(w\in W,\,u\in U)$ with an undirected edge $\{w,u\},$ and every undirected edge $\{u,w\}$ with an arc $\overrightarrow{uw}$ (see Figure~\ref{fig05}).}
{ \begin{thm}\label{thm5.01}
The mixed graph $M_{G'}$ obtained from $M_G$ by the two-way switching is cospectral with $M_G.$
\end{thm}
\begin{proof}
We use a similarity transformation with the diagonal matrix $D$ whose $(v,v)$-entry $D_v$ is equal to $1$ if $v\in W$ and $\bar{\omega}$ if $v\in U$. Let $N=N(M_G)$. The entries of the matrix $N':=D^{-1}ND$ are given by
$$
N'_{uv}=D_u^{-1}N_{uv}D_v.
$$
It is clear that $N'$ is the $N$-matrix for $M_{G'}$. As $N'$ is similar to $N(M_G)$, $M_{G'}$ is cospectral with $M_G.$
\end{proof}
There is a more general switching preserving the spectrum and conserving the underlying graph, based on the structure in Theorem \ref{thm4.1}(i).}

Suppose that the vertex set of $M_G$ is partitioned into six (possibly empty) sets,
\[\label{eq:5.01}
  V(M_G)=V_1\cup V_{-1}\cup V_{\omega}\cup V_{\bar{\omega}}\cup V_{-\omega}\cup V_{-\bar{\omega}}.
\]
An arc $\overrightarrow{xy}$ or an undirected edge $\{x,y\}$ is said to be of \textit{type} $(j,k)$ for $j,k\in \mathbb{T}_6$ if $x\in V_j$ and $y\in V_k$. The partition is said to be \textit{admissible} if both of the following two conditions hold:
\begin{wst}
\item[{\rm (i)}] each undirected edge is one of the type $(j,j),\,(j,\omega\cdot j)$ for $j\in\mathbb{T}_6;$
\item[{\rm (ii)}] each arc is one of the type $(j,j),\,(j,\bar{\omega}\cdot j)$ or $(j,-\omega\cdot j)$ for $j\in\mathbb{T}_6.$
\end{wst}

\begin{figure}[h!]
\begin{center}
\includegraphics[width=100mm]{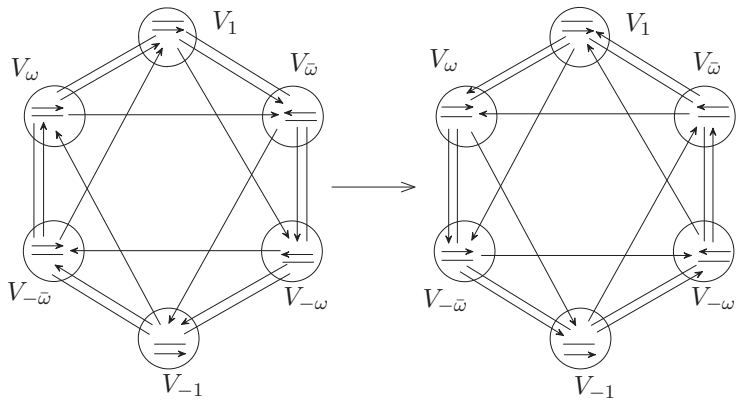} \\
  \caption{Three-way switching with respect to an admissible partition.}\label{fig4}
\end{center}
\end{figure}

A \textit{three-way switching} with respect to the admissible partition \eqref{eq:5.01} is the operation of changing $M_G$ into the mixed graph $M_{G'}$ by making the changes in what follows (see Figure~\ref{fig4}):
\begin{wst}
\item[{\rm (i)}] replacing each undirected edge of type $\left(j,\omega\cdot j\right)$ with an arc directed from $V_j$ to $V_{\omega\cdot j}$ for $j\in\mathbb{T}_6;$
\item[{\rm (ii)}] replacing each arc of type $(j,\bar{\omega}\cdot j)$ with an undirected edge for $j\in\mathbb{T}_6;$
\item[{\rm (iii)}] reversing the direction of each arc of type $(j,-\omega\cdot j)$ for $j\in\mathbb{T}_6.$
\end{wst}

\begin{remark}{\rm
Note that the two-way switching is a special case of the three-way switching in which four of the six sets of the partition \eqref{eq:5.01} are empty, but the three-way switching is not always a sequence of two-way switchings. For example, we can switch the directed hexagon (all the edges of the hexagon are arcs with the same direction) to its underlying graph by a three-way switching, but we cannot switch the directed hexagon to its underlying graph by a sequence of two-way switchings.
}
\end{remark}

Given a mixed graph $M_G$, let $M_G^c$ be its { \textit{converse}} (the mixed graph obtained by reversing all the arcs of $M_G$). It is immediate from the definition of the $N$-matrix that $N(M_G^c)=N(M_G)^T$. This implies the following result.
{ \begin{lem}\label{lem5.2}
A mixed graph $M_G$ and its converse are cospectral.
\end{lem}}
Two mixed graphs $M_{G_1}$ and $M_{G_2}$ are called \textit{switching equivalent} if one can be obtained from the other by a sequence of three-way switchings and operations of taking the converse. A mixed graph $M_G$ is \textit{determined by its spectrum} if $M_G$ is switching equivalent to each mixed graph to which it is cospectral.

Our next result characterizes the mixed graph being cospectral to its underlying graph, the line follows \cite[Theorem 4.1]{0007}.
\begin{thm}\label{thm5.3}
Let $G$ be a connected simple graph of order $n$ and let $M_{G_1}$ be a mixed graph whose underlying graph $G_1$ is a spanning subgraph of $G$. Then the following statements are equivalent:
\begin{wst}
\item[{\rm (a)}]$G$ and $M_{G_1}$ are cospectral.
\item[{\rm (b)}]$\lambda_1(G)=\lambda_1(M_{G_1}).$
\item[{\rm (c)}]$G_1=G$, and the vertex set of $M_{G_1}$ has a partition $\bigcup_{j\in \mathbb{T}_6}V_j$ such that the following holds: for $j\in \mathbb{T}_6$, the induced subgraph $M_{G_1}[V_j]$ contains only undirected edges; each of the rest edges $uv$ of $M_{G_1}$ is an arc $\overrightarrow{uv}$ with $u\in V_j$ and $v\in V_{\bar{\omega}\cdot j}$ for some $j\in \mathbb{T}_6$.
\item[{\rm (d)}]$G$ and $M_{G_1}$ are switching equivalent.
\end{wst}
\end{thm}
\begin{proof}
Clearly, (a) implies (b), and (d) implies (a). By the definition of three-way switching, (c) implies (d) directly. Hence, it suffices to show that (b) implies (c).

{Assume that (b) holds. Let $N=N(G),\, N'=N(M_{G_1})$ and let
$\mathbf{y}=(y_1,\,\ldots,\,y_n)^T\in \mathbb{C}^n$  be a normalized eigenvector of $N'$ corresponding to $\lambda_1(M_{G_1}).$ Following the same line as the proof of \cite[Theorem 4.1]{0007} gives
\[ \label{eq:5.3}
\lambda_1(M_{G_1})\leq \lambda_1(G).
\]
The equality in \eqref{eq:5.3} holding implies $G_1=G$ and
\[\label{eq:5.4}
\, (N')_{jk}\overline{y_j}y_k=|(N')_{jk}||y_j||y_k|
\]
for every edge $v_jv_k$. Since $\mathbf{y}\neq0$, without loss of generality, we can assume that $y_1\in\mathbb{R}^+$, one has $y_1/|y_1|=1$. Then in view of Eq. (\ref{eq:5.4}), we can see, if $v_k\in N_{M_{G_1}}^0(v_1),$ then $N_{1k}'=1,\, y_k/|y_k|=1$; if $v_k\in N_{M_{G_1}}^+(v_1)$, then $N_{1k}'=\omega,\, y_k/|y_k|=\bar{\omega}$; if $v_k\in N_{M_{G_1}}^-(v_1),$ then $N_{1k}'=\bar{\omega},\, y_k/|y_k|=\omega.$

Note that $G_1$ is connected. Then repeating the above argument shows that $y_k/|y_k|\in\mathbb{T}_6$ for $k\in\{1,\ldots, n\}$. Let $V_j=\{v_k\in V(M_G):y_k/|y_k|=j\},\, j\in\mathbb{T}_6.$ Then they construct a partition of $V(M_G)$. It is straightforward to check that the edges within and between the parts are as claimed in (c).
}
\end{proof}
\section{\normalsize Characterizing mixed graphs with rank $2$ or $3$}\setcounter{equation}{0}
When we say the $H$-rank of a mixed graph, we mean the rank of its Hermitian adjacency matrix, and when we say the rank of a mixed graph, we mean the rank of its $N$-matrix.

{Mohar \cite{0007} determined all the mixed graphs with $H$-rank $2$, and constructed a class of mixed graphs which can not be determined by their Hermitian spectra; Wang et al. \cite{0012} considered this problem on the mixed graphs with $H$-rank $3$. For the adjacency rank of signed directed graphs and complex unit gain graphs, one may be referred to \cite{YL2021,PW}.} Inspired directly from \cite{0007,0012}, we are to characterize all the mixed graphs with rank $2$ and $3$, respectively. Furthermore, we show that each connected mixed graph with rank $2$ {is switching equivalent to each connected mixed graph to which it is cospectral.} However, this does not hold for some connected mixed graphs with rank $3$.

Let $M_G$ be a mixed graph of order $n$, the rank of $M_G$ is denoted by $Rank(M_G)$, and the nullity of the $N$-matrix for $M_G$ is denoted by $\eta(M_G)$. Then it is clear that $\eta(M_G)=n-Rank(M_G)$. Thus we can use nullity instead of rank in some cases.

It is well known that $\eta(T)=n-2\mu(T)$ for any tree $T$ of order $n$, where $\mu(T)$ is the matching number of $T$. {According to Corollary \ref{cor2.1},} for any mixed forest, its spectrum is the same as the adjacency spectrum of its underlying graph. {It is a well-known result in gain graph theory (see \cite{NR})}. So we immediately get the following two lemmas, which are the same as the corresponding results for Hermitian adjacency matrices for mixed graphs \cite{0012}.
\begin{lem}\label{lem6.1}
If $M_T$ is a mixed tree of order $n$, then $\eta(M_T)=n-2\mu(T)$, where $\mu(T)$ is the matching number of $T$.
\end{lem}
\begin{lem}\label{lem6.2}
Let $M_P$ be a mixed path of order $n$. Then
\[\label{eq:6.1}
\, \eta(M_P)=\left\{
                         \begin{array}{ll}
                           1, & \text{if}\ n\ \text{is\ odd}, \\
                           0, & \text{if}\ n\ \text{is\ even}.
                         \end{array}
                       \right.
\]
\end{lem}
{ Similar to \cite[Lemma 3.3]{0012}, we have the following result.}
\begin{lem}\label{lem6.3}
Let $M_G$ be a mixed graph containing a pendant edge $uv$, and let $M_{G'}=M_G-u-v$. Then
$
\eta(M_G)=\eta(M_{G'}).
$
\end{lem}
\begin{proof}
Without loss of generality, assume that $v$ is a pendant vertex of $M_G,$ then the lemma can be proved by considering the ranks of $N(M_G)$ and $N(M_G-E),$ where $E$ is the set of edges between $V(M_{G'})$ and $u$ in $M_G$.
\end{proof}
{With the aid of Theorem~\ref{thm3.2}, the nullity of the $N$-matrix for a mixed cycle can be identified.}
\begin{lem}\label{lem6.4}
Let $M_C$ be a mixed cycle of order $n$. Then
\[\label{eq:6.2}
\, \eta(M_C)=\left\{
                         \begin{array}{ll}
                           0, & \text{if $n$ is odd}, \\
                           2, & \text{if $n \equiv2 \pmod{4}$ and $M_C$ is negative}, \\
                           0, & \text{if $n \equiv2 \pmod{4}$ and $M_C$ is positive, semi-positive or semi-negative}, \\
                           2, & \text{if $n \equiv0 \pmod{4}$ and $M_C$ is positive}, \\
                           0, & \text{if $n \equiv0 \pmod{4}$ and $M_C$ is negative, semi-positive or semi-negative}.
                         \end{array}
                       \right.
\]
\end{lem}
\begin{proof}
Denote the characteristic polynomial of $N(M_C)$ by $P_{M_C}(x)=\sum_{j=0}^nc_jx^{n-j}$. {To prove $\eta(M_C)=0$, it is sufficient to prove that $c_n\neq0$; whereas to prove $\eta(M_C)=2$, it is sufficient to prove that $c_{n-2}\neq0$ and $c_{n-1}=c_n=0$}.

{If $n$ is odd, then the only elementary subgraph of $M_C$ with $n$ vertices is $M_C$ itself, and so by Theorem~\ref{thm3.2},
$$
c_n=(-1)^{-1+l_{sn}(M_C)+l_n(M_C)}\cdot 2^{l_p(M_C)+l_n(M_C)}\neq0.
$$

If $n \equiv2 \pmod{4}$ and $M_C$ is negative, then $M_C$ has three $n$-vertex elementary subgraphs, one of which is $M_C$ itself and two of which consist of $\frac{n}{2}$ (oriented) edges, respectively; no $(n-1)$-vertex elementary subgraph; $\frac{n^2}{4}$ elementary subgraphs of order $(n-2),$ each of which consists of $\frac{n-2}{2}$ (oriented) edges. And so by Theorem~\ref{thm3.2},
$$
c_n=(-1)^{-1+1}\cdot 2^1+2\cdot(-1)^{-n+n-\frac{n}{2}}\cdot2^0=0;\,\,\,
c_{n-1}=0;\,\,\,
c_{n-2}=(-1)^{-n+n-\frac{n-2}{2}}\cdot\frac{n^2}{4}\cdot 2^0=\frac{n^2}{4}.
$$

The check of the rest cases is left to the readers.}
\end{proof}
The following lemma is similar to \cite[Lemma 5.1]{0007} {and it can be proved by applying interlacing theorem (see Corollary \ref{cor2.5}).}
\begin{lem}\label{lem6.5}
Suppose that $M_G$ is a mixed graph and $M_{G'}$ is an induced mixed subgraph of $M_G$. Then $Rank(M_G)\ge Rank(M_{G'})$.
\end{lem}
\subsection{\normalsize Mixed graphs with rank $2$}
\begin{lem}\label{lem6.6}
Suppose that $M_G$ is a mixed graph with rank $2$. Then $M_G$ has the following properties:
\begin{wst}
\item[{\rm (a)}]$M_G$ consists of one connected component with more than one vertex together with some isolated vertices.
\item[{\rm (b)}]Every induced subgraph of $M_G$ has rank $0$ or $2$.
\item[{\rm (c)}]The underlying graph of $M_G$ contains no {vertex} induced path on at least $4$ vertices and no {vertex} induced cycle of length at least $5$.
\end{wst}
\end{lem}
The proof of this lemma is similar to that of \cite[Lemma 5.2]{0007}, { so} we omit it here.
\begin{thm}\label{thm6.7}
If $M_G$ is a connected mixed graph with rank $2$, then $G$ is a complete bipartite graph.
\end{thm}
\begin{proof}
According to Lemmas~\ref{lem6.4} and~\ref{lem6.5}, we know that $G$ contains no odd cycle, hence $G$ is bipartite. A shortest path between any two nonadjacent vertices in opposite parts of the bipartition would induce a path on at least $4$ vertices. Since $G$ has no induced $P_4$ (based on Lemma~\ref{lem6.6}), there are no such nonadjacent vertices. Since $G$ contains at least one edge, it is necessarily a complete bipartite graph.
\end{proof}
{Let $u,\,v$ be in $V(M_G)$. Then $u$ and $v$ are \textit{twins} if $M_G$ is switching equivalent to a mixed graph $M'_G$ satisfying $N'_{uj}=N'_{vj}$ for all $j$, where $N'=N(M_G')$. Let $T_{M_G}$ be obtained from $M_G$ by removing all but one from every set of twins.} The following observation is easy to obtain, and enables us to assume that there are no twins when one classifies mixed graphs of a fixed rank.
\begin{lem}\label{lem6.8}
Let $M'_G$ and $M''_G$ be two mixed graphs with the same underlying graph. Then they are switching equivalent if and only if $T_{M'_G}$ and $T_{M''_G}$ are switching equivalent{ ;} $M'_G$ and $T_{M'_G}$ have the same rank.
\end{lem}
\begin{thm}\label{thm6.9}
Let $M_G$ be a mixed graph of order $n$ whose rank is equal to $2$. Then $M_G$ is switching equivalent to $K_{a,b}\cup tK_1$.
\end{thm}
\begin{proof}

{Without loss of generality, assume that $M_G$ is connected (based on Lemma \ref{lem6.6}(a)). By Theorem~\ref{thm6.7}, $G=K_{a,b}$ for $b\geq a\geq1.$ Now either $a=1,\,G$ is a tree and the claim holds, or $a\geq2.$ In the latter case, note that the collection of quadrangles, which must all have weight $1$ by Lemma~\ref{lem6.4}, form a basis of the cycle space. It follows that the weights of all cycles are $1,$ so $M_G$ is switching equivalent to $G$ (see also \cite{NR}).}
\end{proof}
{ \begin{thm}\label{thm6.10}
Let $M_G$ be a connected mixed graph with rank $2.$ Then $M_G$ is switching equivalent to each connected mixed graph to which it is cospectral.
\end{thm}}
\begin{proof}
Let $M_G$ be a connected mixed graph of order $n$ with rank $2$. Then $M_G$ is switching equivalent to $K_{a,b}\, (a\geq b)$. If there exists a connected mixed graph $M_{G'}$ with the same spectrum to $M_G$, then $M_{G'}$ is switching equivalent to $K_{a',b'}\, (a'\geq b')$. {Now
$$
a+b=n=a'+b',\, \, \, \sqrt{ab}=\rho(M_G)=\rho(M_{G'})=\sqrt{a'b'}
$$
imply $a=a',\,b=b'$, i.e., $M_G$ is switching equivalent to $M_{G'}$.}
\end{proof}
In Theorem~\ref{thm6.10}, if the condition ``connected'' is omitted, then { the result is not true}. For example, $K_{4,9}\cup(n-13)K_1$ is cospectral with $K_{6,6}\cup(n-12)K_1$. Note that if $K_{a,b}\cup(n-a-b)K_1$ is cospectral with $K_{a',b'}\cup(n-a'-b')K_1$, then $K_{ta,sb}\cup(n-ta-sb)K_1$ is cospectral with $K_{ta',sb'}\cup(n-ta'-sb')K_1$ for every integer $t,s\geq1$. This implies the following proposition.
\begin{prop}\label{prop1}
There are infinitely many mixed graphs with rank $2$ which are not determined by their spectrum.
\end{prop}
\subsection{\normalsize Mixed graphs with rank $3$}

\begin{figure}[h!]
\begin{center}
\includegraphics[width=120mm]{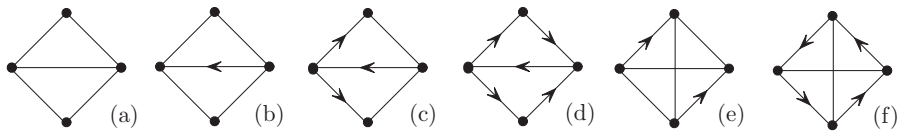} \\
  \caption{Some $4$-vertex mixed graphs.}\label{fig6}
\end{center}
\end{figure}

{ The following result is obtained analogously to Lemma~\ref{lem6.6}.}
\begin{lem}\label{lem6.11}
Suppose that $M_G$ is a mixed graph with rank $3$. Then $M_G$ has the following properties:
\begin{wst}
\item[{\rm (a)}]$M_G$ consists of one connected component with more than one vertex together with some isolated vertices.
\item[{\rm (b)}]Every induced subgraph of $M_G$ has rank $0$, $2$ or $3$.
\item[{\rm (c)}]The underlying graph of $M_G$ contains no {vertex} induced path on at least $4$ vertices and no {vertex} induced cycle of length at least $5$.
\end{wst}
\end{lem}
\begin{lem}\label{lem6.12}
Let $M_G$ be a connected mixed graph of order $4$. Then $Rank(M_G)=3$ if and only if $M_G$ is switching equivalent to one of the mixed graphs as depicted in Figure~\ref{fig6}.
\end{lem}
\begin{proof}
If $M_G$ is a mixed tree or a mixed cycle of order $4$, then $Rank(M_G)=2$ or $4$ by Lemmas~\ref{lem6.1} and ~\ref{lem6.4}. So, we can assume that $M_G$ contains a mixed triangle $M_{C_3}$ with vertex set $\{v_1,v_2,v_3\}$. Let $v$ be in $V(M_G)\setminus \{v_1, v_2, v_3\}$. If $d_G(v)=1,$ then assume that $N_G(v)=\{v_1\}$.  By Lemma~\ref{lem6.3},
\[
  \text{$\eta(M_G)=\eta(M_G-v-v_1)=\eta(K_2)=0,$}
\]
which implies $Rank(M_G)=4$. So we only need to consider {$d_G(v)=2,\,3$.

Denote the characteristic polynomial of $N(M_G)$ by
$$
P_{M_G}(x)=x^4+c_1x^3+c_2x^2+c_3x+c_4.
$$
Then $Rank(M_G)=3$ is equivalent to $c_4=0$ and $c_3\neq 0$. By Lemmas~\ref{lem6.4} and~\ref{lem6.5}, $Rank(M_G)\geq3$, and so it suffices to show $c_4=0$.

In fact, by Theorem~\ref{thm3.2} one has
\[\label{eq:6.06}
c_4=\sum_{M_{G'}}(-1)^{r(G')+l_{sn}(M_{G'})+l_n(M_{G'})}\cdot 2^{l_p(M_{G'})+l_n(M_{G'})},
\]
where the summation is over all spanning elementary subgraphs $M_{G'}$ of $M_G$ and $G'$ is the underlying graph of $M_{G'}$.

If $d_G(v)=2,$ then $G=K_{1,1,2}.$ In this case there are exactly two perfect matchings in $G$ and one spanning mixed cycle (say $M_{C_4}$) in $M_G$. Then
$$
c_4=2+(-1)^{1+l_{sn}(M_{C_4})+l_n(M_{C_4})}\cdot 2^{l_p(M_{C_4})+l_n(M_{C_4})}=0
$$
if and only if $M_{C_4}$ is a positive cycle, i.e., $M_G$ is switching equivalent to (a), (b), (c) or (d) depicted in Figure~\ref{fig6}.

If $d_G(v)=3$, then $G=K_4$. In this case} there are exactly $3$ perfect matchings, say $E_1,\, E_2,\, E_3$, in $G$ and $3$ spanning mixed cycles, say $M_{C_4},\, M'_{C_4},\, M''_{C_4}$, in $M_G$. Let $\mathcal{M}=\{M_{E_1},\, M_{E_2},\, M_{E_3}\}$, $\mathcal{C}=\{M_{C_4},\, M'_{C_4},\, M''_{C_4}\}$. Then \eqref{eq:6.06} gives
\begin{align*}\notag
c_4=&\sum_{M_{G'}\in \mathcal{M}}(-1)^{r(G')+l_{sn}(M_{G'})+l_n(M_{G'})}\cdot 2^{l_p(M_{G'})+l_n(M_{G'})}\\
    &+\sum_{M_{G'}\in \mathcal{C}}(-1)^{r(G')+l_{sn}(M_{G'})+l_n(M_{G'})}\cdot 2^{l_p(M_{G'})+l_n(M_{G'})} \\ \notag
=&\ 3+\sum_{M_{G'}\in \mathcal{C}}(-1)^{r(G')+l_{sn}(M_{G'})+l_n(M_{G'})}\cdot 2^{l_p(M_{G'})+l_n(M_{G'})}.
\end{align*}
Hence, $c_4=0$ if and only if there are three semi-positive mixed cycles or two positive and one semi-negative mixed cycles in $\mathcal{C}$. Note that there does not exist mixed graph $M_{K_4}$ containing three semi-positive spanning cycles. Furthermore, all the mixed graphs $M_{K_4}$ containing two positive and one semi-negative spanning cycles are switching equivalent to (e) or (f) as depicted in Figure~\ref{fig6}.
\end{proof}

\begin{figure}[h!]
\begin{center}
\includegraphics[width=35mm]{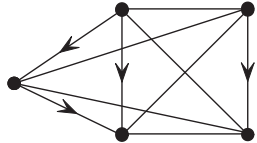} \\
  \caption{The mixed graph $M^0$ whose underlying graph is $K_5$.}\label{fig7}
\end{center}
\end{figure}

\begin{thm}\label{thm6.13}
Let $M_G$ be a connected mixed graph. Then $Rank(M_G)=3$ if and only if $T_{M_G}$ is either a mixed triangle or switching equivalent to $(e)$ or $(f)$ as depicted in Figure~\ref{fig6}.
\end{thm}
\begin{proof}
It suffices to show that any fifth vertex added to one of the mixed graphs in Figure~\ref{fig6} must be a twin. Indeed, note that if vertices $5,\ldots,n$ are local twins of the original four vertices, then adding edges between two twins would strictly increase the rank, it follows that these vertices must also be global twins.

Let $M_G$ be one of the mixed graphs in Figure~\ref{fig6} and suppose $M_{G'}$ is obtained from $M_G$ by adding a vertex $v$ being adjacent to some vertices in $V(M_G)$ such that $Rank(M_{G'})=3.$ Suppose first that $G=K_{1,1,2}.$
\begin{wst}
\item[{\rm(i)}]$d_{G'}(v)=1.$ In this case, there exists a vertex $u$ such that the underlying graph of $M_{G'}-u$ is a triangle with a pendant vertex (henceforth referred to as $\Gamma$), which has exactly one elementary subgraph on $4$ vertices and thus, {by Theorem \ref{thm3.2},} $M_{G'}-u$ has rank $4.$ A contradiction follows by Lemma \ref{lem6.5}.

\item[{\rm(ii)}] $d_{G'}(v)=2$. Then $G'$ contains an induced $\Gamma$ unless two of its vertices (say $y,z$) have degree $4.$ Now, since $M_{G'}-v,\,M_{G'}-u$ and $M_{G'}-w$ ($u,w$ are the non-neighbors of $v$) must all be one of the graphs in Figure~\ref{fig6}(a)-Figure~\ref{fig6}(d), it follows that all triangles in $M_{G'}$ have the same weight, and thus $T_{M_{G'}}$ is a triangle.

\item[{\rm(iii)}]$d_{G'}(v)=3$. Then $G'$ contains an induced $\Gamma$ unless all vertices have degrees at least $3.$ Now every $4$-vertex induced subgraph must be switching equivalent to $M_G,$ and thus $T_{M_{G'}}$ is a triangle.

\item[{\rm(iv)}]$d_{G'}(v)=4.$ This case coincides with (vii), below.
\end{wst}
Next, let $G=K_4.$
\begin{wst}
\item[{\rm(v)}] $d_{G'}(v)=1$ or $d_{G'}(v)=2$. Then $G'$ contains an induced $\Gamma.$

\item[{\rm(vi)}]$d_{G'}(v)=4$. Then $G'=K_5,$ which analogously to above implies that $M_{G'}$ is switching equivalent to $M^0$ as depicted in Figure~\ref{fig7}. Clearly, $M_0$ has rank $5.$

\item[{\rm(vii)}]$d_{G'}(v)=3$. In this case, let $z$ be the non-neighbor of $v.$ Now the graph obtained by removing exactly one neighbor of $v$ must be one of the graphs in Figure~\ref{fig6}(a)-Figure~\ref{fig6}(d), and can be precisely one. This determines the types of two edges incident to $v.$ Similarly, removing a different neighbor of $v$ from $M_{G'}$ determines the final edge, and $v$ is a twin of $z.$
\end{wst}
This completes the proof.
\end{proof}
\begin{remark}{\rm
Wang et al. \cite{0012} characterized all mixed graphs with $H$-rank 3, and showed that {if $M_G$ is a connected mixed graph with $H$-rank $3,$ then $M_G$ is switching equivalent to each connected mixed graph to which it is $H$-cospectral. Here we identify all connected mixed graphs with rank $3$. However, not all connected mixed graphs with rank $3$ have this property.} For example, $K_{8,15,1}$ is not switching equivalent to $M_{K_{3,5,16}}$ whose twin reduction graph is a semi-positive triangle, whereas both of them are cospectral.}
\end{remark}
\section{\normalsize Mixed graphs with small spectral radii}\setcounter{equation}{0}
{Greaves \cite{GG2012} characterized all the gain graphs with gains from the Gauss or Eisenstein integers whose adjacency eigenvalues are contained in $[-2,2]$. The similar problem on $H$-eigenvalues of mixed graphs are considered by Guo and Mohar \cite{0009,0005} and Yuan et al. \cite{BJY}.} Motivated by these nice results, in this section, using interlacing theorem, we characterize all the mixed graphs whose eigenvalues are contained in $(-\alpha,\alpha)$ for $\alpha\in\left\{\sqrt{2},\, \sqrt{3},\, 2\right\}$.

Recall that the spectral radius of an $n$-vertex mixed graph $M_G$ is defined as
$$
\rho(M_G)=\max\{|\lambda_1|, |\lambda_n|\},
$$
where $\lambda_1$ (resp. $\lambda_n$) is the largest (resp. smallest) eigenvalue of $M_G$. Thus the eigenvalues of $M_G$ are contained in $(-\alpha,\alpha)$ if and only if $\rho(M_G)<\alpha$.

\subsection{\normalsize Mixed graphs whose spectral radii are smaller than $\sqrt{3}$}
First we study the case {in which} all eigenvalues are equal to $1$ or $-1$. {This is key for us to determine all the mixed graphs whose spectral radii are smaller than $\sqrt{2}.$} Then we characterize all the mixed graphs whose spectral radii are smaller than $\sqrt{3}.$
\begin{thm}\label{thm7.1}
A mixed graph $M_G$ has the property that $\lambda\in \{-1,1\}$ for each eigenvalue $\lambda$ if and only if $M_G$ is switching equivalent to $tK_2$ for some positive integer $t$.
\end{thm}
The proof of this theorem is the same as \cite[Theorem 9.1]{0005}, which is omitted here.

By Corollary~\ref{cor2.1} we know that all the mixed paths on $n$ vertices are cospectral and all the positive (resp. semi-positive, negative, semi-negative) cycles on $n$ vertices are cospectral. We denote by {$C_n,\, C_n^+,\, C_n^-,\, C_n^=$} the $n$-vertex mixed cycles having no arc, just one arc, just two consecutive arcs with the same direction and just three consecutive arcs with the same direction, respectively. Then they are positive, semi-positive, semi-negative, negative cycles on $n$ vertices, respectively. The following fact is well-known (see { \cite[Section 2.6]{DMS1995}}).
\begin{lem}[\cite{DMS1995}]\label{lem7.2}
The characteristic polynomials of the paths satisfy the recurrence relation $P_{P_n}(x)=xP_{P_{n-1}}(x)-P_{P_{n-2}}(x)$ with $P_{P_0}(x)=1$ and $P_{P_1}(x)=x$. And the spectrum consists of simple eigenvalues
$$
\lambda_j=2\cos\frac{\pi j}{n+1},\,\,\,j=1,\ldots,n.
$$
\end{lem}
{The following result can be obtained by a direct calculation.
\begin{lem}\label{lem7.03}
Let $3\leq n\leq5,$ and let $M_G=M_{C_n}$ be a mixed cycle. Then $\rho(M_G)\geq\sqrt{3},$ unless possibly when $M_G=C_4^=,$ in which case $\rho(M_G)=\sqrt{2}.$
\end{lem}}
We are to characterize all the mixed graphs whose spectral radii are smaller than $\sqrt{2}$ and $\sqrt{3},$ respectively.
\begin{thm}\label{thm7.3}
For a mixed graph {$M_G$,} the followings are equivalent:
\begin{wst}
\item[{\rm (a)}]$\rho(M_G)<\sqrt{2};$
\item[{\rm (b)}]$\rho(M_G)\leq1;$
\item[{\rm (c)}]Every component of $M_G$ is either an undirected edge, an arc or an isolated vertex.
\end{wst}
\end{thm}
\begin{proof}
One may see that (b) implies (a) trivially. Note that, if (c) holds, then together with Theorem~\ref{thm7.1} one has that (b) holds immediately. In order to complete the proof, it suffices to show that (a) implies (c).

{By Lemmas~\ref{lem7.2},~\ref{lem7.03} and Corollary \ref{cor2.5}, $G$ contains no induced $C_3$ or $P_3.$ Hence every 3-vertex induced subgraph of $G$  is $K_2\cup K_1$ or $3K_1$, the conclusion follows.}
\end{proof}
\begin{thm}\label{thm7.4}
Let $M_G$ be an $n$-vertex mixed graph, then $\rho(M_G)<\sqrt{3}$ if and only if every component of $M_G$ is switching equivalent to $P_1,\, P_2,\, P_3,\, P_4$ or $C_4^=$.
\end{thm}
\begin{proof}
``Necessity":\ \ Let $M_G$ be a mixed graph on $n$ vertices with eigenvalues {$\lambda_1,\ldots,\lambda_n$}. Suppose that $\lambda_1<\sqrt{3}$ and $\lambda_n>-\sqrt{3}$. Note that $2\cos{\frac{\pi}{n+1}}$ is increasing as $n$ tends to infinity, and $2\cos{\frac{\pi}{6}}=\sqrt{3}$. Hence, by Corollary~\ref{cor2.5} and Lemma~\ref{lem7.2} we know $M_G$ contains no induced path with order {at least} $5$. As induced mixed cycles with order at least $6$ contain induced paths with order at least $5$, $M_G$ contains no induced mixed cycle with order at least $6$. {By Lemma~\ref{lem7.03},} one may see that $M_G$ contains only $C_4^=$ as an induced mixed cycle {if $M_G$ contains mixed cycles.}

If $M_G$ contains a vertex $v$ with $d_G(v)\geq3$, then $M_G$ contains either an induced mixed star on $4$ vertices or a mixed triangle. Notice that $M_G$ contains no mixed triangle. Hence, $M_G$ must contain an induced mixed star on $4$ vertices. As every mixed star on $4$ vertices is switching equivalent to its underlying graph $K_{1,3}$, and by a direct calculation we know that $\rho(K_{1,3})=\sqrt{3}.$ By Corollary~\ref{cor2.5}, this does not occur for $M_G$. Thus $d_G(v)\leq2$ for all $v\in V(G)$. Therefore, every component of $M_G$ is switching equivalent to $P_1,\, P_2,\, P_3,\, P_4$ or $C_4^=$.

``Sufficiency":\ \ It is straightforward to check that if every component of $M_G$ is switching equivalent to $P_1,\, P_2,\, P_3,$ $P_4$ or $C_4^=$, then $\rho(M_G)<\sqrt{3}$, as desired.
\end{proof}

\subsection{\normalsize Mixed graph whose spectral radii are smaller than $2$}
{In this subsection, we describe all the mixed graphs whose spectral radii are smaller than $2$. Namely, we are to identify all the mixed graphs whose spectra are contained in $(-2,2).$

A $T$-\textit{shape tree} $Y_{a,b,c}$ is a tree with exactly one vertex of degree greater than two such that the removal of this vertex gives rise to paths $P_a, P_b$ and $P_c$.
This tree has $a+b+c+1$ vertices and contains a unique vertex of degree $3$.

The following lemma is well known (see also Lemmens and Seidel \cite{0017} and Smith \cite{0018}).}
\begin{lem}\label{lem7.5}
The largest adjacency eigenvalue of a connected simple graph is smaller than $2$ if and only if the graph is either a path or {the graph $Y_{a,b,1}$,} where either $b=1$ and $a\geq 1$, or $b=2$ and $2\leq a\leq4$.
\end{lem}
As a mixed tree is cospectral with its underlying graph, whose spectral radius is equal to its largest eigenvalue. We have
\begin{cor}\label{cor7.6}
Let $M_G$ be a mixed forest. Then $\rho(M_G)<2$ if and only if each component of $G$ is either a path or {the graph $Y_{a,b,1}$,} where either $b=1$ and $a\geq 1$, or $b=2$ and $2\leq a\leq4$.
\end{cor}
In the following, we consider the case when $M_G$ contains at least one mixed cycle. The spectral radius of $C_n$ is $2$ for $n\geq3$, which follows from the following result.
\begin{lem}[\cite{0003}]\label{lem7.7}
For $n\geq3$, the spectrum of $C_n$ consists of eigenvalues
$$
\lambda_j=2\cos\frac{2j\pi}{n},\,\,\,j=1,\ldots,n.
$$
\end{lem}
\begin{lem}\label{lem7.9}
If $n\geq3$ is odd, then $\lambda$ is an eigenvalue of $C_n$ if and only if $-\lambda$ is an eigenvalue of $C_n^=$.
\end{lem}
\begin{proof}
Let $P_{P_n}(x)=x^n+c_1x^{n-1}+\cdots+c_{n-1}x+c_n,\, P_{P_{n-2}}(x)=x^{n-2}+c'_1x^{n-3}+\cdots+c'_{n-3}x+c'_{n-2}$. As $P_n$ is bipartite, $c_{2j-1}=c'_{2k-1}=0,\, j\in\{1,2,\ldots,\frac{n+1}{2}\},\, k\in\{1,2,\ldots,\frac{n-1}{2}\}.$ Hence $P_{P_n}(x)$ and ${P_{P_{n-2}}(x)}$ are odd functions in $x$. {By Theorem \ref{thm3.7},} one has
$$
P_{C_n}(\lambda)=0\Leftrightarrow P_{P_n}(\lambda)-P_{P_{n-2}}(\lambda)=2\Leftrightarrow P_{P_n}(-\lambda)-P_{P_{n-2}}(-\lambda)=-2\Leftrightarrow P_{C_n^=}(-\lambda)=0,
$$
as desired.
\end{proof}
Together with Lemmas~\ref{lem7.7} and \ref{lem7.9}, we obtain that the spectral radius of $C_n^=$ is $2$ for odd $n$. The following result follows directly from Corollary~\ref{cor2.5}.
\begin{cor}\label{cor7.10}
If $M_G$ is a mixed graph with $\rho(M_G)<2$, then $M_G$ contains no induced positive or odd negative cycle.
\end{cor}
For the other types of mixed cycles, we can also show that their spectral radii are strictly less than $2$, which reads as the following result.
\begin{lem}\label{lem7.11}
If $M_C$ is a semi-positive or semi-negative mixed cycle {of arbitrary order or a} negative mixed cycle of even order, then the spectral radius of $M_C$ is strictly less than $2$.
\end{lem}
\begin{proof}
From \eqref{eq:5.3}, we have $\lambda_1(M_C)\leq\lambda_1(C)$. Similarly,  replacing $\lambda_1(M_C)$ by $|\lambda_n(M_C)|$ in \eqref{eq:5.3} gives  $|\lambda_n(M_C)|\leq\lambda_1(C)$, i.e., $\rho(M_C)\leq\rho(C)=2$. It is sufficient to show that neither $2$ nor $-2$ is an eigenvalue of $M_C$ if $M_C$ is one of those mixed cycles. This follows directly by substituting $2$ and $-2$ into its characteristic polynomial, { which is obtained by a straightforward application of Theorem \ref{thm3.7}.} (Note that if $n$ is even, then $P_{C_n}(2)=P_{C_n}(-2)=0$; if $n$ is odd, then $P_{C_n}(2)=P_{C_n^=}(-2)=0$.)
\end{proof}
{ \begin{lem}\label{lem7.19}
Let $G$ be a connected graph with girth $k,\, k\geq7,$ then $M_G$ has $\rho(M_G)<2$ only if $G= C_k.$
\end{lem}
\begin{proof}
{Suppose that $G\not=C_k$. Then $V(G)\setminus V(C_k)\not= \emptyset$. Choose $u\in V(G)\setminus V(C_k)$ such that $u$ is adjacent to some vertices of $C_k$.

As the girth of $G$ is $k$ $(\geq7)$, one obtains that $u$ is adjacent to just one vertex, say $v$, on $C_k$. Choose a vertex $w$ on $C_k$ such that the distance between $v$ and $w$ on $C_k$ is 3. This gives that $M_G[(V(M_{C_k})\backslash\{w\})\cup\{u\}]$ is switching equivalent to $Y_{k-4,2,1}$. By Corollary~\ref{cor7.6}, $k-4\leq4$, i.e., $k=7$ or~$8$.

If $k=7$, by Corollary~\ref{cor7.10} and Lemma~\ref{lem7.11}, $M_{C_7}$ is semi-positive or semi-negative. By a direct calculation, in both cases, the spectral radius of $M_G[V(M_{C_7})\cup\{u\}]$ is $2.072$, a contradiction.

If $k=8$, then we denote the vertex, say $x$, on $C_8$ satisfying the distance between $v$ and $x$ is 4 on the cycle. This gives $M_G[(V(M_{C_8})\backslash\{x\})\cup\{u\}]$ is switching equivalent to $Y_{3,3,1}$. By Corollary~\ref{cor7.6}, it does not occur.}
\end{proof}
}
\begin{figure}[h!]
\begin{center}
\includegraphics[width=120mm]{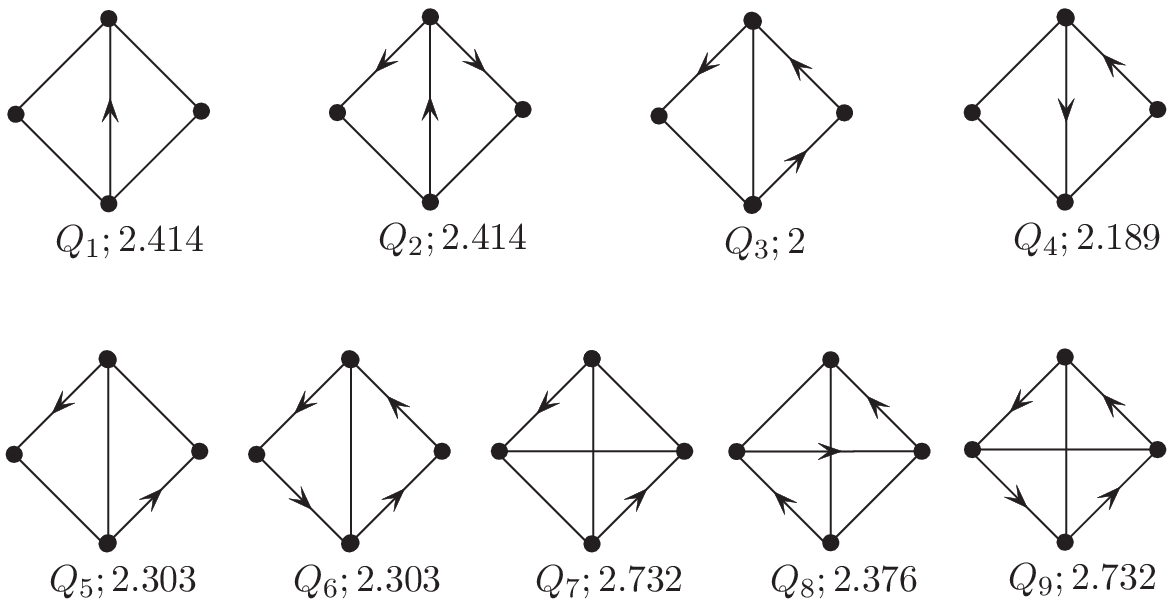} \\
  \caption{Mixed graphs $Q_1,\ldots, Q_9$ together with their spectral radii.}\label{fig8}
\end{center}
\end{figure}
{ \begin{lem}\label{lem7.12}
Let $M_G$ be a mixed graph with $\rho(M_G)<2.$ Then every quadrangle in $M_G$ is chordless and non-positive.
\end{lem}
\begin{proof}
Note that all the 4-vertex graphs $G'$ containing a quadrangle are just $C_4,\, K_{1,1,2}$ and $K_4$. Hence, {by Corollary~\ref{cor7.10}, $M_G$ does not contain chordless positive cycles and negative triangles.}

If $G'=C_4$, then $M_G$ is switching equivalent to one of $C_4^+,\, C_4^-$ and $C_4^=.$ By Lemma \ref{lem7.11}, in any case, $\rho(M_G)<2.$

If $G'=K_{1,1,2},$ then $M_G$ is switching equivalent to one of $Q_i,\,i\in\{1,2,3,4,5,6\}$; see Figure~\ref{fig8}. By a direct calculation, $\rho(Q_1)=\rho(Q_2)=2.414,\,\rho(Q_3)=2,\,\rho(Q_4)=2.189,\,\rho(Q_5)=\rho(Q_6)=2.303.$ All of these give a contradiction by Corollary \ref{cor2.5}.

If $G'=K_4,$ then $M_G$ is switching equivalent to one of $Q_i,\,i\in\{7,8,9\}$; see Figure~\ref{fig8}. By a direct calculation, $\rho(Q_7)=\rho(Q_9)=2.732,\, \rho(Q_8)=2.376.$ All of these also {yield} a contradiction by Corollary \ref{cor2.5}.
\end{proof}
}
\begin{figure}[h!]
\begin{center}
\includegraphics[width=70mm]{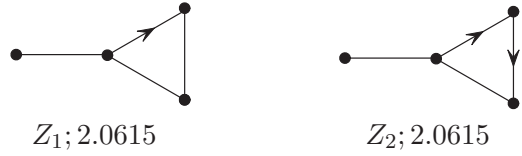} \\
  \caption{Mixed graphs $Z_1$ and $Z_2$ together with their spectral radii.}\label{fig-9}
\end{center}
\end{figure}

{The next lemma characterizes mixed graphs containing a triangle whose spectral radius is less than 2.}
\begin{lem}\label{lem7.13}
Let $M_G$ be a connected mixed graph containing a triangle. Then $\rho(M_G)<2$ if and only if $M_G$ is a semi-positive triangle or a semi-negative triangle.
\end{lem}
\begin{proof}
By Corollary~\ref{cor7.10} and Lemma~\ref{lem7.11}, the triangle contained in $M_G$ is semi-positive or semi-negative. If the order of $M_G$ is $3$, the result is clearly true.

If the order of $M_G$ is $4$, {let $M_{C_3}$ be a triangle in $M_G$,} and let $v$ be a vertex of $M_G$ outside $M_{C_3}$. Then $v$ is adjacent to exactly one vertex in $V(M_{C_3})$, otherwise $M_G$ contains a quadrangle which is not {chordless}. {By Lemma~\ref{lem7.12}, one obtains $\rho(M_G)\geq 2$, a contradiction.} Hence, $M_G$ is switching equivalent to $Z_1$ or $Z_2$; see Figure~\ref{fig-9}. By a direct calculation we obtain that $\rho(Z_1)=\rho(Z_2)=2.0615$, which is a contradiction by Corollary \ref{cor2.5}.
\end{proof}
\begin{lem}\label{lem7.14}
Suppose that $M_G$ is a mixed graph with $\rho(M_G)<2$. Then $\Delta(G)\leq 3$.
\end{lem}
\begin{proof}
Suppose to the contrary that there exists a vertex $v$ in $G$ such that { $d_G(v)\geq 4$.} Consider the mixed graph $Z$ induced by $v$ and four of its neighbors. If $Z$ contains a triangle, then by Lemma~\ref{lem7.13}, $\rho(Z)\geq 2$; if $Z$ contains no triangle, then $Z$ is switching equivalent to a simple bipartite graph $K_{1,4}$. By a direct calculation, $\rho(K_{1,4})=2$. Hence, $\rho(Z)=2$. By Corollary~\ref{cor2.5}, we obtain that $\rho(M_G)\geq \rho(Z)\geq 2$, a contradiction.
\end{proof}

\begin{figure}
\begin{center}
\includegraphics[width=120mm]{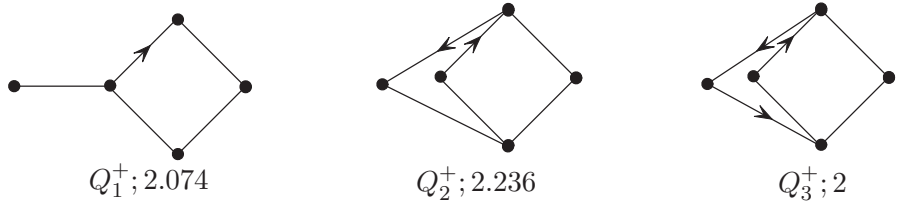} \\
  \caption{Mixed graphs $Q_1^+,\, Q_2^+$ and $Q_3^+$ together with their spectral radii.}\label{fig9}
\end{center}
\end{figure}

\begin{lem}\label{lem7.15}
{Let $M_G$ be a connected mixed graph that contains a semi-positive quadrangle. Then $\rho(M_G)<2$ if and only if $M_G$ is switching equivalent to $C_4^+.$}
\end{lem}
\begin{proof}
{ By Lemma~\ref{lem7.11}, $\rho(C_4^+)<2.$

Let $Q$ be a semi-positive quadrangle contained in $M_G.$} By Lemma \ref{lem7.13}, {$M_G$ does not contain a triangle, and so $Q$ is chordless}. Suppose that $V(M_G)\setminus V(Q)\not= \emptyset.$ Then choose $v\in V(M_G)\setminus V(Q)$ such that $v$ is adjacent to some vertices of $Q$. By Corollary \ref{cor2.5}, one has $\rho(M_G[\{v\}\cup Q])<2.$ As $M_G$ contains no triangle, $v$ does not have three {or more} neighbours in $V(Q)$.

If $v$ has just one neighbor in $V(Q)$, then the mixed graph induced on $V(Q)\cup\{v\}$ is switching equivalent to $Q_1^+$ as depicted in Figure~\ref{fig9}. By a direct calculation, $\rho(Q_1^+)=2.074,$ which is a contradiction.

Next we consider that $v$ has two neighbors in $V(Q)$. The vertex $v$ is adjacent to two non-adjacent vertices of $Q$. Thus, $M_G[\{v\}\cup V(Q)]$ contains three quadrangles. By Corollary~\ref{cor7.10}, $M_G[\{v\}\cup V(Q)]$ contains no positive quadrangle, so there are precisely two distinct switching equivalence classes $Q_2^+$ and $Q_3^+$ displayed in Figure~\ref{fig9}. Direct calculation of their spectral radii leads to a contradiction once more.

Hence $V(M_G)\setminus V(Q)= \emptyset$ and $M_G$ is switching equivalent to $C_4^+.$
\end{proof}
\begin{figure}[h!]
\begin{center}
\includegraphics[width=120mm]{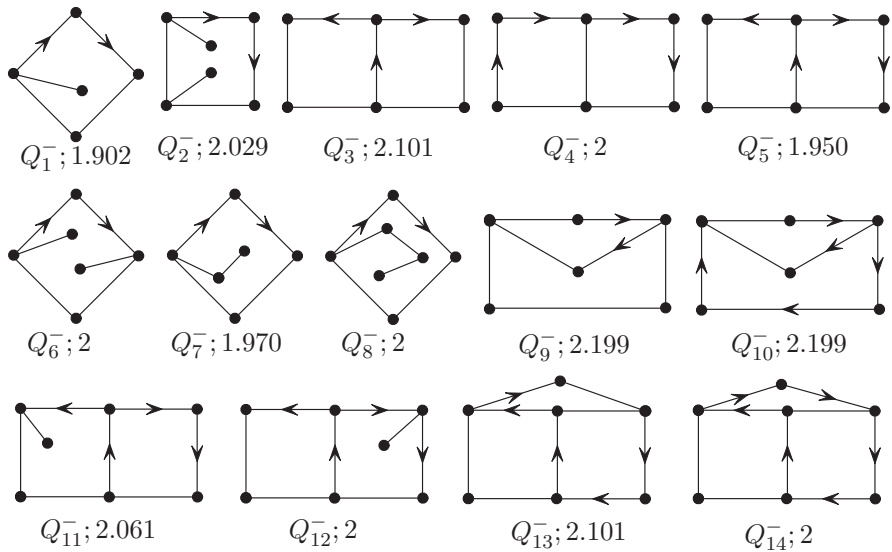} \\
  \caption{Mixed graphs $Q_1^-,\ldots,Q_{14}^-$ together with their spectral radii.}\label{fig-12}
\end{center}
\end{figure}

\begin{lem}\label{lem7.17}
Let $M_G$ be a connected mixed graph containing a semi-negative quadrangle $Q$. If $M_G$ contains no induced subgraph obtained from two semi-negative quadrangles sharing with two consecutive edges. Then $\rho(M_G)<2$ if and only if $M_G$ is switching equivalent to $C_4^-,\, Q_1^-,\, Q_5^-$ or $Q_7^-$, where $Q_1^-,\, Q_5^-$ and $Q_7^-$ are depicted in Figure~\ref{fig-12}.
\end{lem}
{\begin{proof}
By Lemmas~\ref{lem7.11} and~\ref{lem7.12}, $\rho(C_4^-)<2$ and the quadrangle $Q$ is chordless.

Suppose that $V(M_G)\setminus V(Q)\not= \emptyset.$ By Lemmas~\ref{lem7.12}, \ref{lem7.13} and \ref{lem7.15}, $M_G$ contains no triangle, no positive and semi-positive quadrangle. {Note that} $M_G$ contains no other semi-negative quadrangle sharing two consecutive edges with $Q$. {Consequently,} any vertex in $V(M_G)\setminus V(Q)$ has at most one neighbor in $Q.$ Choose $u\in V(M_G)\setminus V(Q)$ such that $u$ is adjacent to a vertex (say $p$) in $V(Q)$. Then the mixed graph induced on $V(Q)\cup\{u\}$ is switching equivalent to $Q_1^-$. By a direct calculation, $\rho(Q_1^-)=1.902$.

Suppose that $V(M_G)\setminus (V(Q)\cup\{u\})\not= \emptyset.$ Then choose $v\in V(M_G)\setminus (V(Q)\cup\{u\})$ such that $v$ is adjacent to some vertices in $V(Q)\cup\{u\}$. Notice that $v$ is adjacent to at most one vertex in $Q$. On the other hand, the maximum degree of $M_G$ is at most $3$ (by Lemma~\ref{lem7.14}), so we have $v\not\sim p$.

If $v$ is adjacent to one vertex in $Q,$ but $v\not\sim u$. Then the mixed graph induced on $V(Q)\cup\{u,v\}$ is switching equivalent to $Q_2^-$ or $Q_6^-$ as depicted in Figure~\ref{fig-12}. By a direct calculation, $\rho(Q_2^-)=2.029,\, \rho(Q_6^-)=2$.

If $v$ is adjacent to $u$, but adjacent to no vertex in $Q$. Then $M_G[V(Q)\cup\{u,v\}]$ is switching equivalent to $Q_7^-$. By a direct calculation, $\rho(Q_7^-)=1.970$.

If $v$ is adjacent to one vertex in $Q,$ and also $v\sim u$. Then $M_G[V(Q)\cup\{u,v\}]$ is switching equivalent to one of $Q_3^-,\,Q_4^-$ and $Q_5^-$ as depicted in Figure~\ref{fig-12} or contains a pentagon. In the former case, by a direct calculation, $\rho(Q_3^-)=2.101,\, \rho(Q_4^-)=2,\, \rho(Q_5^-)=1.950$. In the latter case, by Corollary \ref{cor7.10}, $M_G$ contains no chordless positive and negative pentagon, and so $M_G[V(Q)\cup\{u,v\}]$ is switching equivalent to $Q_9^-$ or $Q_{10}^-$; see Figure~\ref{fig-12}. By a direct calculation we obtain that $\rho(Q_9^-)=\rho(Q_{10}^-)=2.199.$ {From the discussion above, we know that if $V(M_G)\setminus (V(Q)\cup\{u\})\not= \emptyset$ and $v$ is adjacent to some vertices in $V(Q)\cup\{u\}$, then by Corollary \ref{cor2.5} $M_G[V(Q)\cup\{u,v\}]$ must be switching equivalent to $Q_5^-$ or $Q_7^-$.

Suppose that $V(M_G)\setminus (V(Q)\cup\{u,v\})\not= \emptyset.$ Then choose $w\in V(M_G)\setminus (V(Q)\cup\{u,v\})$ such that $w$ is adjacent to some vertices in $V(Q)\cup\{u,v\}$.} First we consider that $M_G[V(Q)\cup\{u,v\}]$ is switching equivalent to $Q_7^-$. As $w$ is adjacent to at most one vertex in $Q$ and by Lemma \ref{lem7.13} $M_G$ contains no triangle, $w$ is adjacent to at most two vertices in $V(Q)\cup\{u,v\}.$ {Together with Corollary \ref{cor7.6} and a similar discussion of $v$ as above}, if $w$ is adjacent to only one vertex in $V(Q)\cup\{u,v\}$, {then} $M_G[V(Q)\cup\{u,v,w\}]$ is switching equivalent to $Q_8^-$ as depicted in Figure~\ref{fig-12}. By a direct calculation, $\rho(Q_8^-)=2$. If $w$ is adjacent to two vertices in $V(Q)\cup\{u,v\},$ then $M_G[V(Q)\cup\{u,v,w\}]$ is switching equivalent to $Q_{12}^-$ as depicted in Figure~\ref{fig-12}. By a direct calculation, $\rho(Q_{12}^-)=2$, which is a contradiction by Corollary \ref{cor2.5}.

Now we consider that $M_G[V(Q)\cup\{u,v\}]$ is switching equivalent to $Q_5^-$. By Lemma \ref{lem7.13} $M_G$ contains no triangle and by Lemma \ref{lem7.14} $M_G$ has maximum degree at most $3,$ $w$ is adjacent to at most two vertices in $V(Q)\cup\{u,v\}.$ If $w$ is adjacent to only one vertex in $V(Q)\cup\{u,v\}$, then $M_G[V(Q)\cup\{u,v,w\}]$ is switching equivalent to $Q_{11}^-$ or $Q_{12}^-$; see Figure~\ref{fig-12}. By a direct calculation, $\rho(Q_{11}^-)=2.061$ and $\rho(Q_{12}^-)=2$. If $w$ is adjacent to two vertices in $V(Q)\cup\{u,v\}$, then by a similar discussion of $v$ as above one has $M_G[V(Q)\cup\{u,v,w\}]$ is switching equivalent to $Q_{13}^-$ or $Q_{14}^-$, see Figure~\ref{fig-12}. By a direct calculation, $\rho(Q_{13}^-)=2.101$ and $\rho(Q_{14}^-)=2$, a contradiction by Corollary \ref{cor2.5}.
\end{proof}}
\begin{figure}[h!]
\begin{center}
\includegraphics[width=120mm]{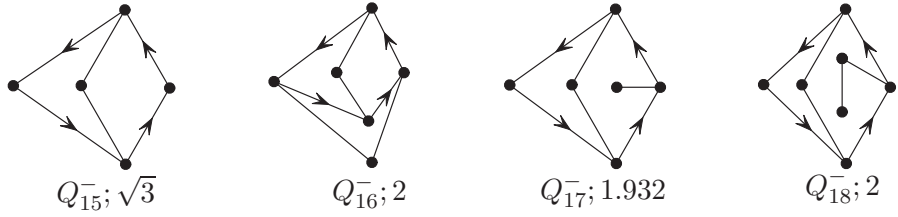} \\
  \caption{Mixed graphs $Q_{15}^-,\, Q_{16}^-,\, Q_{17}^-$, $Q_{18}^-$ together with their spectral radii.}\label{fig12}
\end{center}
\end{figure}
\begin{lem}\label{lem7.18}
Let $M_G$ be a connected mixed graph containing a subgraph obtained from two semi-negative quadrangles sharing with two consecutive edges. Then $\rho(M_G)<2$ if and only if $M_G$ is switching equivalent to $Q_{15}^-$ or $Q_{17}^-$, where $Q_{15}^-$ and $Q_{17}^-$ are depicted in Figure~$\ref{fig12}$.
\end{lem}
\begin{proof}
Note that all mixed graphs obtained from two semi-negative quadrangles sharing with two consecutive edges are switching equivalent to $Q_{15}^-$ {(based on Lemma \ref{lem7.12})}. Without loss of generality, we assume that $Q_{15}^-$ is a mixed subgraph of $M_G.$ In view of Lemma~\ref{lem7.13}, $M_G$ contains no triangle, $Q_{15}^-$ is an induced mixed subgraph of $M_G$. By a direct calculation, $\rho(Q_{15}^-)=\sqrt{3}$.

Suppose that $V(M_G)\setminus V(Q_{15}^-)\not= \emptyset.$ Then choose $u\in V(M_G)\setminus V(Q_{15}^-)$ such that $u$ is adjacent to some vertices of $Q_{15}^-$. { By Lemma \ref{lem7.14},} $u$ is adjacent to at most three vertices in $Q_{15}^-$.

If $u$ is adjacent to two or three vertices in $Q_{15}^-$, then by Lemmas~\ref{lem7.12} and \ref{lem7.15}, $M_G$ contains no positive or semi-positive quadrangle, {which is only possible when} $M_G[V(Q_{15}^-)\cup\{u\}]$ is switching equivalent to $Q_{16}^-$ as depicted in Figure~\ref{fig12}. By a direct calculation, $\rho(Q_{16}^-)=2$, {which is a contradiction by Corollary \ref{cor2.5}}.

If $u$ is adjacent to only one vertex in $Q_{15}^-$, then by Lemma~\ref{lem7.14}, this vertex has degree $2$ in $Q_{15}^-$. $M_G[V(Q_{15}^-)\cup\{u\}]$ is switching equivalent to $Q_{17}^-$ (see Figure~\ref{fig12}). By a direct calculation, $\rho(Q_{17}^-)=1.932$.

If $V(M_G)\setminus (V(Q_{15}^-)\cup\{u\})\not= \emptyset,$ then choose $v\in V(M_G)\setminus (V(Q_{15}^-)\cup\{u\})$ such that $v$ is adjacent to some vertices in $V(Q_{15}^-)\cup\{u\}$. {According to the discussion above, $v$ is adjacent to at most one vertex in $Q_{15}^-,$ and so $v$ is adjacent to at most two vertices in $V(Q_{15}^-)\cup\{u\}$. By Corollary \ref{cor7.10} and Lemma~\ref{lem7.14}, if $v$ is adjacent to two vertices in $V(Q_{15}^-)\cup\{u\},$ then we deduce that $M_G[V(Q_{15}^-)\cup\{u,v\}]$ contains an induced subgraph which is switching equivalent to $Q_9^-$ or $Q_{10}^-$, a contradiction. If $v$ is adjacent to exactly one vertex in $V(Q_{15}^-)\cup\{u\},$ then $M_G[V(Q_{15}^-)\cup\{u,v\}]$ is switching equivalent to $Q_{18}^-$ (see Figure~\ref{fig12}) or contains an induced subgraph which is switching equivalent to $Q_6^-$. By a direct calculation, $\rho(Q_6^-)=\rho(Q_{18}^-)=2$, a contradiction by Corollary \ref{cor2.5}.}
\end{proof}

{In what follows, we are to characterize the mixed graphs $M_G$ containing negative quadrangles with $\rho(M_G)<2.$}

Let $a,\, b,\, c,\, d$ be nonnegative integers. Let $\Box_{a,b,c,d}$ be a mixed graph obtained from a negative quadrangle with consecutive vertices $v_1,\, v_2,\, v_3,\, v_4$ by attaching undirected paths of lengths $a,\, b,\, c,\, d$ to $v_1,\, v_2,\, v_3$ and $v_4$, respectively. This graph has $a+b+c+d+4$ vertices. It is easy to see that the resulting mixed graph is unique up to switching equivalence.

In the discussion of $H$-matrices for mixed graphs, a mixed cycle is negative if and only if its weight is $-1$ (see Guo and Mohar \cite{0009}, Liu and Li \cite{0008}). This definition coincides with ours. By comparing the characteristic polynomials of the $H$-matrix \cite[Theorem 2.8]{0008} and the $N$-matrix (Theorem~\ref{thm3.2}) for a mixed graph, we know that if $M_G$ is a unicyclic mixed graph (i.e, the underlying graph of $M_G$ is a unicyclic graph) with the unique mixed cycle negative, then the characteristic polynomials of theirs are the same. Hence \cite[Lemma 4.11]{0009} gives
\begin{lem}\label{lem7.21}
Let $M_G$ be a unicyclic mixed graph with a negative quadrangle, then $\rho(M_G)<2$ if and only if $M_G$ is switching equivalent to one of the following mixed graphs:
\begin{wst}
\item[{\rm (1)}]$\Box_{a, 0, c, 0}$, where $a\geq c\geq0$;
\item[{\rm (2)}]$\Box_{3, 1, 0, 0},\, \Box_{2, 1, 1, 0},\, \Box_{2, 1, 0, 0},\, \Box_{1, 1, 1, 1},\, \Box_{1, 1, 1, 0}$ or $\Box_{1, 1, 0, 0}$.
\end{wst}
\end{lem}
By a similar discussion as the proof of \cite[Lemma 4.13]{0009}, we obtain the following lemma.
\begin{lem}\label{lem7.22}
Suppose that $M_G$ is a connected mixed graph with $\rho(M_G)<2$, then any two vertices $u$ and $v$ of degree $3$ are at distance at most $3$ in
$G$. 
\end{lem}
{
\begin{lem}\label{lem7.23}
Let $M_G$ be a connected mixed graph that contains at least two chordless mixed cycles, in which at least one is a negative quadrangle. Then $\rho(M_G)<2$ if and only if $M_G$ is switching equivalent to a mixed graph in $\{Q_5^-,\,Q_1^=,\,Q_4^=,\,Q_5^=,\,Q_6^=,\, Q_8^=,\,Q_9^=,\,Q_{10}^=,\,Q_{11}^=\}$ as depicted in Figures~\ref{fig-12} and~\ref{Fig16}.
\end{lem}
\begin{figure}[h!]
\begin{center}
\includegraphics[width=110mm]{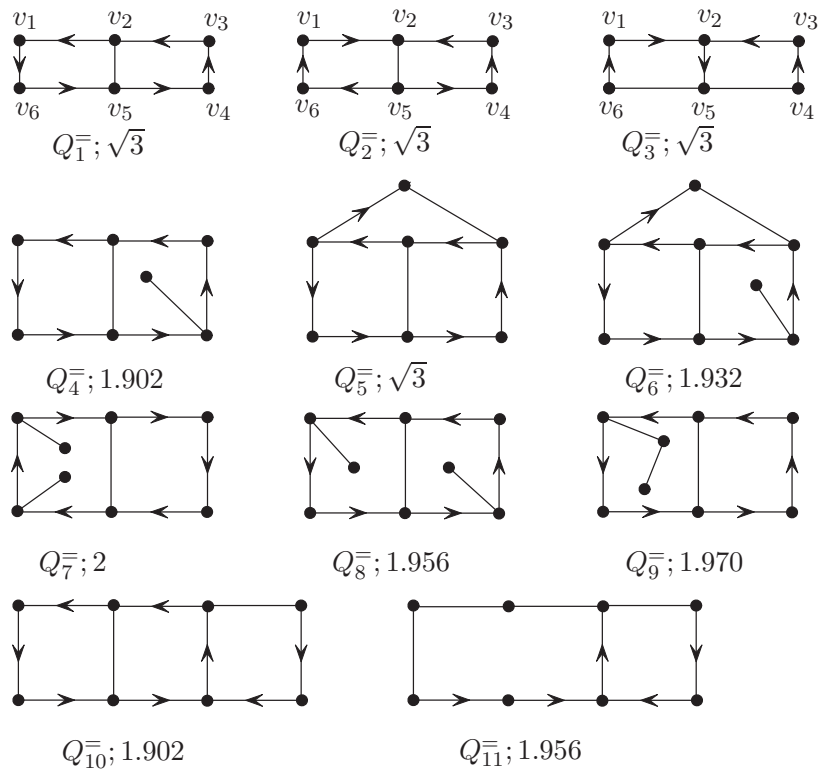} \\
  \caption{Mixed graphs $Q_1^=,\ldots,Q_{11}^=$ together with their spectral radii.}\label{Fig16}
\end{center}
\end{figure}
\begin{proof}
Let $Q_1,\, Q_2$ be two chordless mixed cycles in $M_G$, where $Q_1$ is a negative quadrangle. By Lemmas~\ref{lem7.12}, \ref{lem7.13} and \ref{lem7.15}-\ref{lem7.18}, it is sufficient to consider the cases that $Q_2$ is a negative quadrangle or a mixed cycle of length at least $5.$

Suppose that $Q_1$ and $Q_2$ have no vertex in common. Let $M_P$ be a shortest mixed path between $Q_1$ and $Q_2$ {and denote the end vertices of $M_p$ by $u,v$, where $u\in V(Q_1),\, v\in V(Q_2).$ Take $M_P$ together with two edges in $Q_1$ (resp. $Q_2$) such that these two edges are incident with $u$ (resp. $v$). Clearly, this subgraph is a mixed tree containing two vertices of degree $3$. By Corollary~\ref{cor7.6}, this subgraph has spectral radius at least 2. So it is not an induced subgraph of $M_G$}. By Lemmas~\ref{lem7.13} and \ref{lem7.22}, $M_G$ has no triangle and the length of $M_P$ is at most $3$, and so the only possibility is that a vertex in $Q_1$ is adjacent to a vertex in $Q_2$. Since $M_P$ is a shortest path, it means that $M_P$ is a single edge and thus we have two edges joining adjacent vertices in $Q_1$ with adjacent pair of vertices in $Q_2$. This forms a new quadrangle having common vertices with $Q_2$. So it is enough to consider that $Q_1$ and $Q_2$ have at least one vertex in common.
By Lemma~\ref{lem7.14}, $G$ has maximum degree at most $3$. {Consequently, $Q_1$ and $Q_2$ either have one edge in common or have two consecutive edges in common.

{\bf Case 1.} \ $Q_1$ and $Q_2$ have two consecutive edges in common. }
If $Q_2$ is a negative quadrangle, then $M_G[V(Q_1)\cup V(Q_2)]$ contains three quadrangles. {It is quite simply impossible for all three quadrangles of a $K_{1,3,1}$ to have gains $-1$.}

If $Q_2$ is a mixed cycle of length $k\geq5,$ then $M_G[V(Q_1)\cup V(Q_2)]$ contains two chordless mixed $k$-cycles and one chordless negative quadrangle. Furthermore, one mixed $k$-cycle is positive if and only if the other is negative. In this case, by Corollary~\ref{cor7.10}, the spectral radius of $M_G[V(Q_1)\cup V(Q_2)]$ is at least $2$. { Similarly,} one mixed $k$-cycle is semi-positive if and only if the other is semi-negative. In this case, $M_G[V(Q_1)\cup V(Q_2)]$ is switching equivalent to the mixed graph on the left in Figure~\ref{Fig15}, and the labels at vertices on the right in Figure~\ref{Fig15} show an eigenvector corresponding to the eigenvalue $2$. This implies that the spectral radius of $M_G[V(Q_1)\cup V(Q_2)]$ is at least $2$, a contradiction.
\begin{figure}[h!]
\begin{center}
\includegraphics[width=90mm]{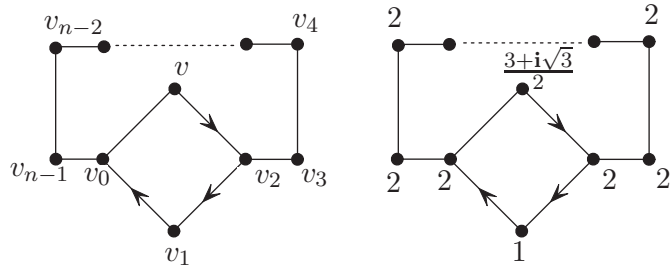} \\
  \caption{Long cycle plus a vertex with two neighbors {at distance two}.}\label{Fig15}
\end{center}
\end{figure}

{{\bf Case 2.} \ $Q_1$ and $Q_2$ have just one edge in common. }

{\bf Subcase 2.1.} \ $Q_2$ is a negative quadrangle. In this subcase, $M_G[V(Q_1)\cup V(Q_2)]$ contains a mixed hexagon $X$ with at least one chord. {Note that $M_G$ contains no triangle, positive quadrangle, semi-positive quadrangle and semi-negative quadrangle. Hence, the endpoints of any chord of $X$ are at distance 3 on $X.$ Clearly, there are at most three such chords of $X.$}

Furthermore, each quadrangle contained in $M_G[V(Q_1)\cup V(Q_2)]$ is negative.
And so if $X$ contains two or three chords, then $X$ is a directed hexagon (all the edges of it are arcs with the same direction) with two or three undirected chords. By a direct calculation, in these two cases, $\rho(M_G[V(Q_1)\cup V(Q_2)])=2$, a contradiction. It remains to consider the case that there is only one chord of $X$.

If the chord is undirected, then $M_G[V(Q_1)\cup V(Q_2)]$ is switching equivalent to $Q_1^=$ or $Q_2^=$, whereas if the chord is directed, then $M_G[V(Q_1)\cup V(Q_2)]$ is switching equivalent to $Q_3^=$, where $Q_1^=,\, Q_2^=$ and $Q_3^=$ are depicted in Figure~\ref{Fig16}.

We are to show that both $Q_2^=$ and $Q_3^=$ are  switching equivalent to $Q_1^=.$ 
In fact, for $Q_3^=$, we can take $U=\{v_5\},\, W=V(Q_3^=)\backslash\{v_5\}$. Then a two-way switching with respect to this partition gives $Q_2^=$.
For $Q_2^=$, we can take $V_1=\{v_2,v_3,v_4,v_5\},\, V_{-\omega}=\{v_6\},\, V_{-\bar{\omega}}=\{v_1\}$ and $V_{-1}=V_{\omega}=V_{\bar{\omega}}=\emptyset.$ Note that this partition is admissible. Hence, a three-way switching with respect to this partition yields $Q_1^=$. {Indeed, the fact that both $Q_2^=$ and $Q_3^=$ are switching equivalent to $Q_1^=$ can be also deduced by considering a basis of the cycle space.} By a direct calculation, one may easily obtain $\rho(Q_1^=)=\sqrt{3}.$

In what follows, we are to characterize all the mixed graphs $M_G$ containing $Q_1^=$ as an induced subgraph with $\rho(M_G)<2.$

{
Suppose that $V(M_G)\setminus (V(Q_1)\cup V(Q_2))\not= \emptyset$. Choose $u\in V(M_G)\setminus (V(Q_1)\cup V(Q_2))$ such that $u$ is adjacent to some vertices in $V(Q_1)\cup V(Q_2)$. By Lemmas~\ref{lem7.13} and~\ref{lem7.14}, $u$ is adjacent to at most $2$ vertices in $V(Q_1)\cup V(Q_2)$. Consequently, we consider the following two subcases.
\begin{wst}
\item[{\rm(i)}] $u$ is adjacent to exactly one vertex in $V(Q_1)\cup V(Q_2)$. In this subcase, $M_G[V(Q_1)\cup V(Q_2)\cup\{u\}]$ is switching equivalent to $Q_4^=$ (see Figure~\ref{Fig16}). By a direct calculation, $\rho(Q_4^=)=1.902.$

\item[{\rm(ii)}] $u$ is adjacent to two vertices in $V(Q_1)\cup V(Q_2)$. In this subcase, $M_G[V(Q_1)\cup V(Q_2)\cup\{u\}]$ either contains a chordless mixed pentagon sharing two consecutive edges with a negative quadrangle and which has been discussed in Case 1, or is switching equivalent to $Q_5^=$. By a direct calculation, $\rho(Q_5^{=})=\sqrt{3}.$
\end{wst}

From the discussion above, we know that if $V(M_G)\setminus (V(Q_1)\cup V(Q_2))\not= \emptyset$ and $u$ is adjacent to some vertices in $V(Q_1)\cup V(Q_2)$, then $M_G[V(Q_1)\cup V(Q_2)\cup\{u\}]$ must be switching equivalent to $Q_4^=$ or $Q_5^=$. Suppose that $V(M_G)\setminus (V(Q_1)\cup V(Q_2)\cup\{u\})\not= \emptyset.$ Then choose $v\in V(M_G)\setminus (V(Q_1)\cup V(Q_2)\cup\{u\})$ such that $v$ is adjacent to some vertices in $V(Q_1)\cup V(Q_2)\cup\{u\}$.

Bear in mind that $M_G[V(Q_1)\cup V(Q_2)\cup\{u\}]$ is switching equivalent to $Q_4^=$ or $Q_5^=$. By Lemma~\ref{lem7.14}, $v$ is adjacent to at most $3$ vertices in $V(Q_1)\cup V(Q_2)\cup\{u\}$. At first we consider the former. 
\begin{wst}
\item[{\rm(iii)}] $v$ is adjacent to exactly one vertex in $V(Q_1)\cup V(Q_2)\cup\{u\}.$ Then $M_G[V(Q_1)\cup V(Q_2)\cup\{u,v\}]$ either contains an induced mixed tree with two vertices of degree $3$, by Corollary~\ref{cor7.6}, this {does not} happen; or is switching equivalent to one of $Q_7^=,\,Q_8^=,\,Q_9^=$ (see Figure~\ref{Fig16}). By a direct calculation, $\rho(Q_7^=)=2,\, \rho(Q_8^=)=1.956$ and $\rho(Q_9^=)=1.970$.

\item[{\rm(iv)}] $v$ is adjacent to two vertices in $V(Q_1)\cup V(Q_2)$. Then this case is the same as the case {\rm(ii)}.

\item[{\rm(v)}] $v\sim u$ and $v$ is adjacent to exactly one vertex in $V(Q_1)\cup V(Q_2)$. Then $M_G[V(Q_1)\cup V(Q_2)\cup\{u,v\}]$ either contains a chordless mixed pentagon or hexagon and which will be discussed in Subcases 2.2 and 2.3; or is switching equivalent to $Q_{10}^=$ (see Figure~\ref{Fig16}). By a direct calculation, $\rho(Q_{10}^=)=1.902$.
\end{wst}

Now we consider that $M_G[V(Q_1)\cup V(Q_2)\cup\{u\}]$ is switching equivalent to $Q_5^=$.  
\begin{wst}
\item[{\rm(vi)}] $v$ is adjacent to exactly one vertex in $V(Q_1)\cup V(Q_2)\cup\{u\}.$ In this subcase, by Lemma~\ref{lem7.14}, one has $M_G[V(Q_1)\cup V(Q_2)\cup\{u,v\}]$ is switching equivalent to $Q_6^=$ (see Figure~\ref{Fig16}). By a direct calculation, $\rho(Q_6^=)=1.932$.

\item[{\rm(vii)}] $v$ is adjacent to two vertices in $V(Q_1)\cup V(Q_2)\cup\{u\}.$ In this subcase, a direct calculation gives that the spectral radius of $M_G[V(Q_1)\cup V(Q_2)\cup\{u,v\}]$ is $2$, a contradiction.

\item[{\rm(viii)}] $v$ is adjacent to three vertices in $V(Q_1)\cup V(Q_2)\cup\{u\}.$ In this subcase, not all quadrangles of $M_G[V(Q_1)\cup V(Q_2)\cup\{u,v\}]$ are negative, a contradiction.
\end{wst}

Repeat the procedures as above, {one may see that each graph $M_G$ among $\{Q_6^=,\, Q_8^=,\, Q_9^=,\, Q_{10}^=\}$ is maximal with respect to $\rho(M_G)< 2$}. This completes the proof of the case that $Q_2$ is a negative quadrangle.

{\bf Subcase 2.2.} \ $Q_2$ is a mixed pentagon. Then by Corollary \ref{cor7.10}, $Q_2$ is a semi-positive or semi-negative pentagon. Note that the case in which $M_G$ contains a pentagon sharing two consecutive edges with a negative quadrangle has been discussed in Case 1. Hence, it suffices to consider the case $E(M_G[V(Q_1)\cup V(Q_2)])=E(Q_1)\cup E(Q_2).$ Let $u$ be a vertex in $V(Q_1)\backslash V(Q_2),$ then by a direct calculation,  one may obtain that $M_G[V(Q_2)\cup\{u\}]$ has spectral radius $2.076$ whenever $Q_2$ is either semi-positive or semi-negative, a contradiction.

{\bf Subcase 2.3.} \ $Q_2$ is a mixed hexagon. Then by Corollary \ref{cor7.10}, $Q_2$ is a semi-positive, semi-negative or negative hexagon. By a similar discussion as that in Subcase 2.2, it suffices to consider $E(M_G[V(Q_1)\cup V(Q_2)])=E(Q_1)\cup E(Q_2).$

Let $u$ be a vertex in $V(Q_1)\backslash V(Q_2).$ Then we obtain $M_G[V(Q_2)\cup\{u\}]$ has spectral radius $2.037$ if $Q_2$ is semi-positive and $M_G[V(Q_2)\cup\{u\}]$ has spectral radius $2$ if $Q_2$ is semi-negative, a contradiction.

Now we consider $Q_2$ is negative. Then, $M_G[V(Q_1)\cup V(Q_2)]$ is switching equivalent to $Q_{11}^=$. By a direct calculation, $\rho(Q_{11}^=)=1.956$.
It is straightforward to check that $Q_{11}^=$ is maximal with respect to the spectral radius being less than 2.

{\bf Subcase 2.4.} \ $Q_2$ is a mixed cycle with length at least $7$. Then by a similar discussion as that in subcase 2.2, it is sufficient for us to consider the case $E(M_G[V(Q_1)\cup V(Q_2)])=E(Q_1)\cup E(Q_2).$ Let $u$ be in $V(Q_1)\backslash V(Q_2).$ By Lemma \ref{lem7.19}, $M_G[V(Q_2)\cup\{u\}]$ has spectral radius at least $2$, a contradiction.}

By Cases 1 and 2, we complete the proof.
\end{proof}

{In the next two lemmas we characterize the mixed graphs $M_G$ of girth 5 or 6 satisfying $\rho(M_G)<2,$ respectively.}
\begin{lem}\label{lem7.16}
Let $M_G$ be a mixed graph, where $G$ is a connected graph with girth $5$. Then $\rho(M_G)<2$ if and only if $M_G$ is a semi-positive pentagon or a semi-negative pentagon.
\end{lem}
\begin{proof}
As the girth of $G$ is $5$, the pentagons contained in $M_G$ are chordless. By Corollary~\ref{cor7.10} and Lemma~\ref{lem7.11}, the pentagons contained in $M_G$ are semi-positive or semi-negative. 

If $V(M_G)=5$, the result is clear true. If $V(M_G)\ge 6$,, let $M_{C_5}$ be a mixed pentagon contained in $M_G$, and let $v$ be a vertex of $M_G$ outside $M_{C_5}$. As $M_G$ contains no triangle, no quadrangle, $v$ is adjacent to exactly one vertex in $V(M_{C_5})$. A direct calculation gives $\rho(M_G)=2.076$ (when the mixed pentagons contained in $M_G$ are either semi-positive or semi-negative), {which is a contradiction by Corollary \ref{cor2.5}}.
\end{proof}}

\begin{figure}[h!]
\begin{center}
\includegraphics[width=90mm]{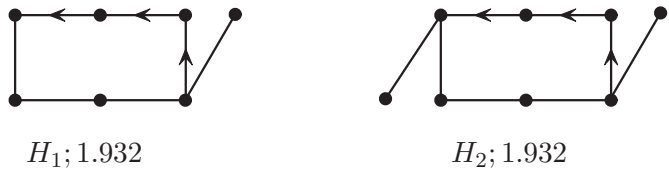} \\
  \caption{Mixed graphs $H_1,\,H_2$ together with their spectral radii.}\label{Fig13}
\end{center}
\end{figure}
{ \begin{lem}\label{lem7.20}
Let $M_G$ be a mixed graph, where $G$ is a connected graph with girth $6$. Then $\rho(M_G)<2$ if and only if $M_G$ is switching equivalent to $C_6^+,\, C_6^-,\, C_6^=,\, H_1$ or $H_2$, where $H_1$ and $H_2$ are depicted in Figure~\ref{Fig13}.
\end{lem}
\begin{proof}
Notice that the girth of $G$ is $6$. Hence the hexagons contained in $G$ are chordless. By Corollary~\ref{cor7.10} and Lemma~\ref{lem7.11}, the mixed hexagons contained in $M_G$ are semi-positive or semi-negative or negative. 
Let $Z$ be a mixed hexagon contained in $M_G$.

If $V(M_G)\setminus V(Z)= \emptyset$, the result is clearly true. Suppose now that $V(M_G)\setminus V(Z)\not= \emptyset.$ Then choose $u\in V(M_G)\setminus V(Z)$ such that $u$ is adjacent to some vertices of $Z$. As the girth of $G$ is $6,\, u$ is adjacent to just one vertex of $Z$. If $Z$ is a semi-positive (resp. semi-negative) hexagon, then by a direct calculation, the spectral radius of $M_G[V(Z)\cup\{u\}]$ is $2.074$ (resp. $2$), a contradiction. 

So we consider that $Z$ is a negative hexagon. 
As $u$ is adjacent to exact one vertex, say $z_1,$ of $Z,$ then $M_G[V(Z)\cup\{u\}]$ is switching equivalent to $H_1$ as depicted in Figure~\ref{Fig13}. By a direct calculation, $\rho(H_1)=1.932$.

Suppose that $V(M_G)\setminus (V(Z)\cup\{u\})\not= \emptyset.$ Choose $v\in V(M_G)\setminus (V(Z)\cup\{u\})$ such that $v$ is adjacent to some vertices in $V(Z)\cup\{u\}.$ By a similar discussion as $u,$ one obtains that $v$ is adjacent to exactly one vertex of $Z$, and so $v$ is adjacent to at most two vertices in $V(Z)\cup\{u\}.$

If $v$ is adjacent to exactly one vertex, say $z_2,$ in $V(Z)\cup\{u\}$, then by Corollary \ref{cor7.6} and Lemma \ref{lem7.14}, $z_2$ is on $Z,$ and $z_1,z_2$ are at distance 3 on $Z$. Hence, $M_G[V(Z)\cup\{u,v\}]$ is switching equivalent to $H_2$; see Figure~\ref{Fig13}. By a direct calculation, $\rho(H_2)=1.932$.

If $v$ is adjacent to two vertices in $V(Z)\cup\{u\}$. As the girth of $G$ is $6,\, v$ is adjacent to $u$ and $z_2,$ where $z_2$ and $z_1$ are at distance 3 on $Z.$ Then $M_G[V(Z)\cup\{u,v\}]$ contains three chordless hexagons, and not all of them are negative. {By the discussion above this does not occur.}

If $V(M_G)\setminus (V(Z)\cup\{u,v\})\not= \emptyset,$ then $M_G[V(Z)\cup\{u,v\}]$ is switching equivalent to $H_2$. Choose $w\in V(M_G)\setminus (V(Z)\cup\{u,v\})$ such that $w$ is adjacent to some vertices in $V(Z)\cup\{u,v\}.$ With the same discussion as $v$, we deduce a contradiction in this case.
\end{proof}

Together with above lemmas, the following conclusion is clear.}
\begin{thm}\label{thm7.24}
Let $M_G$ be a connected mixed graph. Then $\rho(M_G)<2$ if and only if $M_G$ is switching equivalent to one of the following:
\begin{wst}
\item[{\rm (a)}]$C_n^+,\, C_n^-;$
\item[{\rm (b)}]$C_n^=$ with even $n;$
\item[{\rm (c)}]$P_n;$
\item[{\rm (d)}]$Y_{a,b,1},$ where either $b=1$ and $a\geq 1,$ or $b=2$ and $2\leq a\leq4;$
\item[{\rm (e)}]$Q_1^-,\, Q_5^-,\, Q_7^-;$ see Figure~$\ref{fig-12};$
\item[{\rm (f)}]$Q_{15}^-,\,Q_{17}^-;$ see Figure~$\ref{fig12};$
\item[{\rm (g)}]$\Box_{a, 0, c, 0},$ where $a\geq c\geq0;$
\item[{\rm (h)}]$\Box_{3, 1, 0, 0},\, \Box_{2, 1, 1, 0},\, \Box_{2, 1, 0, 0},\, \Box_{1, 1, 1, 1},\, \Box_{1, 1, 1, 0},\, \Box_{1, 1, 0, 0};$
\item[{\rm (i)}]$Q_1^=,\, Q_4^=,\,Q_5^=,\,Q_6^=,\, Q_8^=,\,Q_9^=,\,Q_{10}^=,\,Q_{11}^=;$ see Figure~$\ref{Fig16};$
\item[{\rm (j)}]$H_1,\, H_2;$ see Figure~$\ref{Fig13}.$
\end{wst}
\end{thm}

\begin{remark}
{\rm From Theorem~\ref{thm7.4}, one may see that there are only finitely many connected mixed graphs with all eigenvalues in the interval $(-\sqrt{3},\, \sqrt{3})$. But from Theorem~\ref{thm7.24}, we can see that there are infinite many connected mixed graphs with all eigenvalues in the interval $(-2,\, 2)$.}
\end{remark}

\section*{\normalsize Acknowledgment}
We take this opportunity to thank the anonymous reviewers for their critical reading of the manuscript and suggestions which have immensely helped us in getting the article to its present form. The authors would like to thank one of the anonymous reviewers for the constructive suggestions and inputs for the arrangement of results in Section 6 and the proofs of Theorems~\ref{thm6.9}, \ref{thm6.13} and Lemma~\ref{lem7.12}.


\begin{thebibliography}{99}
\small \setlength{\itemsep}{-.8mm}
\bibitem{AAH2016}A. Abiad, W.H. Haemers, Switched symplectic graphs and their $2$-ranks, Des. Codes Cryptogr. 81 (1) (2016) 35-41.
\bibitem{0016} C. Adiga, B.R. Rakshith, W. So, On the mixed adjacency matrix of a mixed graph, Linear Algebra Appl. 495 (2016) 223-241.
\bibitem{ABC2012}B. Arsi\'c, D. Cvetkovi\'c, S.K. Simi\'c, M. \v{S}kari\'c, Graph spectral techniques in computer sciences, Appl. Anal. Discrete Math. 6(1) (2012) 1-30.
\bibitem{0003} R.B. Bapat, Graphs and Matrices, Springer, New York, 2010.
\bibitem{NB1993}N. Biggs, Algebraic Graph Theory, Cambridge Mathematical Library, 2nd edn., Cambridge University Press, Cambridge, 1993.
\bibitem{AEB2012}A.E. Brouwer, W.H. Haemers, Spectra of Graphs, Universitext, Springer, New York, 2012.
\bibitem{BR2010}R.A. Brualdi, Spectra of digraphs, Linear Algebra Appl. 432 (2010) 2181-2213.
\bibitem{MC2010}M. Cavers, S.M. Cioab\v{a}, S. Fallat, D.A. Gregory, W.H. Haemers, S.J.Kirkland, J.J. McDonald, M. Tsatsomeros, Skew-adjacency matrices of graphs, Linear Algebra Appl 436(12) (2012) 4512-4529.
\bibitem{CSM2017}S.M. Cioab\v{a}, W.H. Haemers, J.R. Vermette, The graphs with all but two eigenvalues equal to $-2$ or 0, Des. Codes Cryptogr. 84 (1-2) (2017) 153-163.
\bibitem{DMS1995}D. Cvetkovi\'c, M. Doob, H. Sachs, Spectra of Graphs, Theory and Applications, 3rd edn., Johann Ambrosius Barth, Heidelberg, 1995.
\bibitem{0002} D. Cvetkovi\'c, P. Rowlinson, S. Simi\'c, An Introduction to the Theory of Graph Spectra, New York: Cambridge University Press; 2009.
\bibitem{CDS2011}D. Cvetkovi\'c, S.K. Simi\'c, Graph spectra in computer science, Linear Algebra Appl. 434 (6) (2011) 1545-1562.{
\bibitem{ALG}A. L. Gavrilyuk, S. Suda, On the multiplicities of digraph eigenvalues, arXiv: 1911.11055 (2019).}
\bibitem{GGR2001}C. Godsil, G. Royle, Algebraic Graph Theory, vol. 207 of Graduate Texts in Mathematics, Springer-Verlag, New York, 2001.
{ \bibitem{GG2012}G. Greaves, Cyclotomic matrices over the Eisenstein and Gaussian integers, J. Algebra 372 (2012) 560-583.}
\bibitem{0009} K. Guo, B.J. Mohar, Digraphs with Hermitian spectral radius below $2$ and their cospectrality with paths, Discrete Math. 340 (11) (2017) 2616-2632.
\bibitem{0005} K. Guo, B.J. Mohar, Hermitian adjacency matrix of digraphs and mixed graphs, J. Graph Theory  85 (1) (2017) 217-248.
\bibitem{P4}Hs.H. G\"{u}nthard, H. Primas, Zusammenhang von Graphentheorie und MO-Theorie von Molekeln mit Systemen konjugierter Bindungen, Helv. Chim. Acta 39 (1956) 1645-1653.
\bibitem{0017} P.W.H. Lemmens, J.J. Seidel, Equiangular lines, J. Algebra 24 (1973) 494-512.
\bibitem{LS}S.C. Li, W.T. Sun, On split graphs with three or four distinct (normalized) Laplacian eigenvalues, J. Combin. Des. 28 (11) (2020) 763-782.
\bibitem{SCL2020}S.C. Li, W. Wei, The multiplicity of an $A_\alpha$-eigenvalue: A unified approach for
mixed graphs and complex unit gain graphs, Discrete Math. 343 (2020) 111916.
\bibitem{LSG2012}X.L. Li, Y.T. Shi, I. Gutman, Graph Energy, Springer, New York, 2012.
\bibitem{0008} J.X. Liu, X.L. Li, Hermitian-adjacency matrices and Hermitian energies of mixed graphs, Linear Algebra Appl. 466 (2015) 182-207.
{ \bibitem{YL2021}Y. Lu, J. Wu, Bounds for the rank of a complex unit gain graph in terms of its maximum degree, Linear Algebra Appl. 610 (2021) 73-85.}
\bibitem{0001} B.J. Mohar, A new kind of Hermitian matrices for digraphs, Linear Algebra Appl. 584 (2020) 343-352.
\bibitem{0007} B.J. Mohar, Hermitian adjacency spectrum and switching equivalence of mixed graphs, Linear Algebra Appl.  489  (2016) 324-340.
\bibitem{0004} J. Oxley, D. Vertigan, G. Whittle, On maximum-sized near-regular and $\sqrt[6]{1}$-matroids, Graphs Combin. 14 (2) (1998) 163-179.
{ \bibitem{NR} N. Reff, Spectral properties of complex unit gain graphs, Linear Algebra Appl. 436 (9) (2012) 3165-3176.
\bibitem{NR2016} N. Reff, Oriented gain graphs, line graphs and eigenvalues, Linear Algebra Appl. 506 (2016) 316-328.
\bibitem{AS2019} A. Samanta, M.R. Kannan, On the spectrum of complex unit gain graph, arXiv: 1908.10668 (2019).}
\bibitem{0018} J.H. Smith, Some properties of the spectrum of a graph, in: Combinatorial Structures and their Applications (Proc. Calgary Internat. Conf., Calgary. Alta., 1969), Gordon and Breach, New York, 1970, pp. 403-406.
\bibitem{ZS2015}Z. Stani\'c, Inequalities for Graph Eigenvalues, London Mathematical Society Lecture Note Series, 423. Cambridge University Press, Cambridge, 2015. xi+298 pp.
\bibitem{0012} Y. Wang, B.J. Yuan, S.D. Li, C.J. Wang, Mixed graphs with $H$-rank $3$, Linear Algebra Appl. 524 (2017) 22-34.
\bibitem{DW1996}D.B. West, Introduction to Graph Theory, Prentice Hall, Inc., Upper Saddle River, NJ, 1996.
\bibitem{0006} G. Whittle, On matroids representable over $GF(3)$ and other fields, Trans. Amer. Math. Soc. 349 (2) (1997) 579-603.
{ \bibitem{PW2020} P. Wissing, E.R. van Dam, The negative tetrahedron and the first infinite family of connected digraphs that are strongly determined by the Hermitian spectrum, J. Combin. Theory Ser. A, 173 (2020).
\bibitem{PW} P. Wissing, E.R. van Dam, Spectral fundamentals and characterizations of signed directed graphs, arXiv: 2009.12181 (2020).
\bibitem{BJY} B.J. Yuan, Y. Wang, S.C. Gong, Y. Qiao, On mixed graphs whose Hermitian spectral radii are at most $2$, Graphs Combin. 36 (5) (2020) 1573-1584.}
\end{thebibliography}
\end{document}